\documentclass{article}

\usepackage[french,british]{babel}

\usepackage[T1]{fontenc}

\usepackage{amsmath}

\usepackage{amsfonts}
\usepackage{amsbsy}
\usepackage{amssymb}

\usepackage{amsthm}

\theoremstyle{plain}

\usepackage{latexsym}
\usepackage{amsmath}
\usepackage{amsfonts}
\usepackage{amssymb}
\usepackage{psfrag}
\usepackage{graphicx}

\usepackage{amsmath,amssymb}
\usepackage{graphicx}

\newtheorem{remark}{Remark}
\newtheorem{example}{Example}
\newtheorem{lemma}{Lemma}
\newtheorem{theorem}{Theorem}

\newtheorem{corollary}{Corollary}

\newtheorem{definition}{Definition}

\theoremstyle{Remark}

\theoremstyle{Proposition}
\newtheorem{Proposition}{Proposition}[section]

\theoremstyle{definition}

\newcommand{\lp}[1]{\left( \begin{array}{#1} }
\newcommand{\rp}{\end{array} \right)}

\newcommand{\be}{\begin{equation}}
\newcommand{\ee}{\end{equation}}

\begin{document}

\title{Interpolation problem for periodically correlated stochastic sequences with missing observations}

\author
{Iryna Golichenko\thanks
{Department of Mathematical Analysis and Probability Theory, National Technical  University of Ukraine ''Igor Sikorsky Kyiv Politechnic Institute'', Ukraine,
idubovetska@gmail.com },
Mikhail Moklyachuk\thanks
{Department of Probability Theory, Statistics and Actuarial
Mathematics, Taras Shevchenko National University of Kyiv, Kyiv 01601, Ukraine, Moklyachuk@gmail.com},
}

\date{\today}

\maketitle

\renewcommand{\abstractname}{Abstract}
 \begin{abstract}
The problem of mean square optimal estimation of linear functionals
which depend on the
unobserved values of a periodically correlated stochastic sequence is considered.
The estimates are based on observations of the sequence with a noise.
Formulas for calculation the mean square errors and the spectral
characteristics of the optimal estimates of functionals
are derived in the case of spectral certainty, where the spectral densities of the sequences are
exactly known.
Formulas that determine the least favorable spectral
densities and the minimax spectral characteristics are proposed in the case of spectral uncertainty, where the spectral densities of the sequences are
not exactly known while some classes of admissible spectral densities are specified.
\end{abstract}

\vspace{2ex}

\maketitle

\vspace{2ex}

\section{Introduction}

W.R.~Bennett \cite{Bennet} in 1958  introduced cyclostationarity as a phenomenon describing signals in channels of communication.
Studying the statistical characteristics of information transmission, he calls the group of telegraph signals a cyclostationary process, that is the process whose group of statistics changes periodically with time.
W.A.~Gardner and L.~E.~Franks \cite{Gardner1975}  highlighted the similarity of  cyclostationary processes, which form a subclass of nonstationary processes, with stationary processes.
W.A.~Gardner \cite{Gardner1994}, W.~A.~Gardner, A.~Napolitano and L.~Paura \cite{Gardner2006} presented bibliography of works in which properties and applications of cyclostationary processes were studied. Recent developments and applications of cyclostationary signal analysis are reviewed in the papers by A.~Napolitano  \cite{Nap1}, \cite{Nap2}.  Note that in different sources cyclostationary processes are called periodically stationary, periodically nonstationary, periodically correlated. We will use the term periodically correlated processes.

E.G.~Gladyshev \cite{Glad1961} was the first who analysed the spectral properties and representations of periodically correlated sequences based on its connection with the vector valued stationary sequences. He formulated the necessary and sufficient conditions for determining the periodically correlated sequence in terms of the correlation function. A.~Makagon  \cite{Makagon1999}, \cite{Makagon2011} presented a detailed spectral analysis of periodically correlated sequences. The main ideas of the research  of periodically correlated sequences are outlined in the book by H.~L.~Hurd and A.~Miamee~\cite{Hurd}.

The problem of estimation of unobserved values of random processes is one of the very important and topical subsections of the theory of stochastic processes. Processes that are observed can be completely defined by its characteristics (correlation function, spectral density, canonical decomposition) or their characteristics can be defined only by the set of admissible values of characteristics.
The linear extrapolation and interpolation problems for stationary
stochastic sequences under the condition that the spectral densities are
 exactly known were first investigated by A.~N.~Kolmogorov \cite{Kolmogorov}.
Methods of solutions of the extrapolation and filtering problems for stationary
processes and sequences with rational spectral densities were
developed by N.~Wiener \cite{Wiener} and A.~M.~Yaglom \cite{Yaglom-a, Yaglom-b}. Estimation
problems for vector-valued stationary processes were investigated by Yu.~A.~Rozanov~\cite{Rozanov} and E.~J.~Hannan~\cite{Hannan1974}.

 The basic techniques of statistics of stochastic processes are summarized in the books by  D.~Z.~Arov and H.~Dym~\cite{Arov2018}, L.~Aggoun  and R.~J.~Elliott~\cite{Aggoun2004},
 I.~V.~Basawa  and  B.~L.~S.~Prakasa Rao \cite{Basawa1980},  S.~Cohen and   R.~J.~Elliott~\cite{Cohen2015},
 M.~S.~Grewal and A.~P.~Andrews \cite{Grewal2015}, G.~Kallianpur \cite{Kallianpur1980}, Yu.~A.~Kutoyants  \cite{Kutoyants1998}, \cite{Kutoyants2004},
 R.~S.~Liptser  and A.~N.~Shiryaev  \cite{Liptser2001a,Liptser2001b},
 B.~L.~S.~Prakasa Rao \cite{PrakasaRao2010,PrakasaRao2012}, B.~L.~Rozovsky  and S.~V.~Lototsky~\cite{Rozovsky2018},
 M.~B.~Rajarshi~\cite{Rajarshi2012},
W.~A.~Woodward, H.~L.~Gray and A.~C.~Elliott  \cite{Woodward2017}.
The estimation problems occur in different studies. We refer to D.~V.~Koroliouk \cite{Koroliouk1}, D.~V.~Koroliouk et al.~\cite{Koroliouk2}, D.~V.~Koroliouk  and V.~S.~Koroliuk~\cite{Koroliouk3},  where there is investigated the difference stochastic equation
$ \Delta\alpha_{t+1}=-V_0\alpha_t+\sigma_0 \Delta W_{t+1},\, t\geq0,$
which determines a sequence $\alpha_t$, $t\geq0$, for the stochastic component $\Delta W_{t+1}$, $t\geq0$, and studied the problem of filtration of stationary Gaussian statistical experiments considered for the solution $\alpha_t$, $t\geq0$, of the indicated equation.

Since processes often accompanied by undesirable noise it is naturally to assume that the exact value of  spectral density is unknown and the model of process  is given by a set of restrictions on spectral density.
K~S.~Vastola and H~V.~Poor~\cite{VastPoor}  showed for   certain classes of spectral densities that the Wiener filter is very sensitive to minor changes of spectral model unlike the robust Wiener filter.  That is the filter is the least sensitive to the worst case of uncertainty. Thus,  it is reasonable to use the minimax (robust) estimation method, which allows to define the optimal estimate for all densities from a certain given class of the admissible spectral densities simultaneously.
 Ulf Grenander \cite{Grenander} was the first who
proposed the minimax approach to the extrapolation problem for stationary processes.
A survey of results in minimax-robust methods of data processing can be found in the paper by S.~A.~Kassam and H.~V.~Poor \cite{Kassam}.
Formulation and investigation of the problems of extrapolation, interpolation and filtering of linear functionals which depend on the unknown values of stationary sequences and processes  from observations with and without noise are presented by M.~P.~Moklyachuk in the papers \cite{Moklyachuk2000}--\cite{Moklyachuk2015}.
Similar problems of the optimal estimation for the vector-valued stationary sequences and processes were examined by M.~P.~Moklyachuk \cite{Moklyachuk1993}--\cite{Moklyachuk:1995ROSE} and by M.~P.~Moklyachuk and O.~Yu.~Masyutka  \cite{MoklyachukMas2006a}--\cite{MoklyachukMas2012}.
In their papers M.~M.~Luz and M.~P.~Moklyachuk \cite{luz1}--\cite{luz11} investigated the minimax estimation problems for linear functionals which depends on unobserved values of stochastic sequences with stationary increments.
 P.~S.~Kozak and M.~P.~Moklyachuk~\cite{Kozak} study estimates of functionals constructed from random sequences with periodically stationary increments.
In their papers I.~I.~Golichenko(Dubovets'ka) and M.~P.~Moklyachuk \cite{Dubovetska0}--\cite{Dubovetska9}, \cite{MoklyachukD2016} presented results of investigation of the
interpolation, extrapolation and filtering problems for linear functionals from periodically correlated stochastic sequences and processes.

The prediction problem for stationary sequences with missing observations is investigated in the papers by P.~Bondon \cite{Bondon1,Bondon2},
R.~Cheng, A.~G.~Miamee and M.~Pourahmadi \cite{Cheng1}, R.~Cheng and M.~Pourahmadi~\cite{Cheng2},
Y.~Kasahara, M.~Pourahmadi and A.~Inoue~\cite{Kasahara,Pourahmadi}.
The detailed analysis of the estimation problems with missing observations are presented in the
paper by B.~Abraham\cite{Abraham},
books by
M.~J.~Daniels and J.~W.~Hogan \cite{Daniels2008},
 R.~J.~A.~Little and D.~B.~Rubin \cite{Little2002},
 P.~E.~McKnight et al \cite{McKnight2007}, M.~M.~Pelagatti~\cite{Pelagatti}.

In the papers by M.~P.~Moklyachuk and M.~I.~Sidei
\cite{MSidei2015}--\cite{MSidei2017}
results of investigations of the interpolation, extrapolation and filtering problems for stationary stochastic sequences and processes with missing observations are proposed.
The results of the study of the extrapolation, interpolation and filtering problems for linear functionals constructed from
unobserved values of multidimensional stochastic sequences and processes are presented in the papers by O.~Yu.~Masyutka, M.~P.~Moklyachuk and M.~I.~Sidei
\cite{Masyutka:Moklyachuk:Sidei1}--\cite{Masyutka:Moklyachuk:Sidei5}, \cite{Moklyachuk:Mas:Sidei}.
We also refer to the book by M.~P.~Moklyachuk, O.~Yu.~Masyutka and I.~I.~Golichenko \cite{Moklyachuk:Mas:Gol2018}
where results of the investigation of the problem of mean square optimal estimation (forecasting, interpolation, and filtering) of linear functionals constructed from unobserved values
of periodically correlated isotropic random fields are described.

In this paper we deal with the problem of optimal linear estimation of
the functional $A_s{\zeta}$ which depends on the
unobserved values of a periodically correlated stochastic sequence
${\zeta}(j)$. Estimates are based on observations of the sequence
${\zeta}(j)+{\theta}(j)$ at points  $j\in{\mathbb
Z}\setminus S$,  where $S=\bigcup _{l=0}^{s-1}\{M_l+1,\dots,M_l+N_{l+1}\}$. ${\theta}(j)$ is an uncorrelated
with ${\zeta}(j)$ periodically correlated stochastic sequence.
Formulas for calculation the mean square errors and the spectral
characteristics of the optimal estimates of the functional
$A_{s}\zeta$ are proposed in the case of spectral certainty where the spectral densities are
exactly known. Formulas that determine the least favorable spectral
densities and minimax spectral characteristics are proposed in the case of spectral uncertainty where the spectral densities are
not exactly known while some classes of admissible spectral densities are given.

The paper is organized as follows.
The spectral properties of periodically correlated stochastic sequences and their correlation functions are described in Section 2.
Relations of periodically correlated stochastic sequences with multidimensional stationary sequences are discussed in this section.

In section 3 we consider the problem of mean square optimal linear estimation of the
the functional $$A_{s}\vec {\xi}=\sum_{l=0}^{s-1}\sum_{j=M_l+1}^{M_l+N_{l+1}}\vec {a}^{\top}(j)\vec
{\xi}(j),\,\,M_l=\sum_{k=0}^l(N_k+K_k),\,\,N_0=K_0=0,$$ which depends on the unknown values of a $T$-dimensional
stationary stochastic sequence  $\vec {\xi}(j)$,
based on observations of the sequence $\vec
{\xi}(j)+\vec{\eta}(j)$ at points $j\in {\mathbb Z}
\setminus S$, where $S=\bigcup _{l=0}^{s-1}\{M_l+1,\dots,M_l+N_{l+1}\}$.
Formulas for calculation the mean square error and the spectral
characteristic of the optimal estimate of the functional
$A_{s}\vec \xi$ are proposed in the case where spectral density matrices of the sequences
$\vec {\xi}(j)$ and $\vec{\eta}(j)$ are exactly known.

In section 4 we consider the problem of mean square optimal linear estimation of the
the functional
 $$A_s{\zeta}=\sum_{l=0}^{s-1}\sum_{j=M_l+1}^{M_l+N_{l+1}}{a}(j){\zeta}(j),
\,\,
 M_l=\sum_{k=0}^l(N_k+K_k),\,\,N_0=K_0=0,$$
which depends on the unknown values of T-PC stochastic
sequence $\zeta(j)$, based on
observations of the sequence $\zeta(j)+\theta(j)$ at points $j\in {\mathbb
Z} \setminus S,$ where $S=\bigcup _{l=0}^{s-1}\{M_l+1,\dots,M_l+N_{l+1}\}$.

In section 5 we consider the problem
of optimal estimation for the linear functional
$$A_s{\zeta}=\sum_{l=0}^{s-1}\sum_{j=M_l+1}^{M_l+N_{l+1}}{a}(j){\zeta}(j),\,\,M_l=\sum_{k=0}^l(N_k+K_k),\,\,N_0=K_0=0,$$
which depends on the
unknown values of $T$-PC sequence
${\zeta}(j)$ from observations of the sequence
${\zeta}(j)+{\theta}(j)$ at points $j\in{\mathbb
Z}\setminus S$,
where the number of missed observations at each of the intervals is a multiple of the period $T$.
In sections 4 and 5 the estimation problem is investigated in the case of spectral certainty, where the spectral densities of observed sequences are exactly known.

In section 6 we describe the minimax approach to the problem of  estimation of the linear functionals. In this case we find the
estimate which minimizes the mean square error for all spectral
densities from the given set of admissible densities simultaneously.

In section 7
the least favorable spectral densities and the minimax (robust) spectral characteristics of the optimal estimate of ${A_s \vec \zeta}$ are found  for the class $D_{0}^{-}$ of admissible spectral densities.

In section 8
the least favorable spectral densities and the minimax (robust) spectral characteristics of the optimal estimate of ${A_s \vec \zeta}$ are found  for the class $D_{G}^{-}$ of admissible spectral densities.

\section{Periodically correlated and multidimensional stationary  sequences}

The term \textit {periodically correlated} process was introduced by E.~G.~Gladyshev \cite{Glad1961} while W.~R.~Bennett \cite{Bennet} called random and periodic processes
\textit {cyclostationary} process.

Periodically correlated sequences are stochastic sequences that have
periodic structure (see the book by H.~L.~Hurd and A.~Miamee~\cite{Hurd}).

\begin{definition} A complex valued stochastic
sequence ${\zeta}(n),n\in{\mathbb Z}$ with zero mean,
$\textsf{E}{\zeta}(n)=0$, and finite variance, $\textsf{E}|{\zeta}(n)|^{2}<+\infty$,
is called cyclostationary or periodically
correlated (PC) with period $T$ ($T$-PC) if for every
$n,m\in{\mathbb Z}$
\begin{equation}\label{2.1}
\textsf{E}{\zeta}(n+T)\overline{\zeta(m+T)}=R(n+T,m+T)=R(n,m)
\end{equation}
and there are no smaller values of $T>0$ for
which (\ref{2.1})  holds true.
\end{definition}

\begin{definition} {A complex valued
T-variate stochastic  sequence $\vec{\xi}(n)=\left\{ \xi_{\nu}(n)
\right\}_{\nu = 1 }^{T},$ $n\in {\mathbb Z}$ with zero mean,
$\textsf{E}\xi_{\nu}(n)=0, \nu=1,\dots,T$, and $\textsf{E}||\vec {\xi}(n)||^{2}<\infty$
is called stationary if for all $n ,m\in {\mathbb Z}$ and
$\nu,\mu\in \{1,\dots,T\}$}
\[
\textsf{E}\xi_{\nu}(n)\overline{\xi_{\mu}(m)}=R_{\nu\mu}(n,m)=R_{\nu\mu}(n-m).
\]
\noindent If this is the case, we denote $R(n)=\left\{ R_{\nu\mu}(n)
\right\}_{\nu,\mu = 1 }^{T}$ and call it the \textit {covariance
matrix} of T-variate stochastic sequence $\vec{\xi}(n)$.
\end{definition}

\begin{Proposition} {\rm{(E. G. Gladyshev \cite{Glad1961})}}.
A stochastic sequence $\zeta(n)$ is PC with period $T$ if and only if there exists
a $T$-variate stationary sequence $\vec {\xi}(n)=\left\{ \xi_{\nu}(n)
\right\}_{\nu = 1 }^{T}$ such that $\zeta(n)$ has the representation
\begin{equation}\label{2.2}
\zeta(n)=\sum_{\nu=1}^{T}{e}^{2\pi in\nu/T}\xi_{\nu}(n),\,\,n\in {\mathbb
Z}.
\end{equation}
\noindent The sequence $\vec \xi(n)$ is called  \textit {generating sequence} of the sequence $\zeta(n)$.
\end{Proposition}

\begin{Proposition} {\rm{(E. G. Gladyshev \cite{Glad1961})}}. A complex valued stochastic
sequence ${\zeta}(n),n\in{\mathbb Z}$ with zero mean and finite variance is PC with period $T$ if and only if the  $T$-variate blocked sequence $\vec \zeta(n)$ of the form
\begin{equation} \label{block}
[\vec{\zeta}(n)]_{p}=\zeta(nT+p),\,\,n\in {\mathbb
Z},p=1,\dots,T
\end{equation}
 is stationary.
 \end{Proposition}

We will denote by $f^{\vec \zeta}(\lambda)=\left\{ f^{\vec \zeta}_{\nu\mu}(\lambda)
\right\}_{\nu,\mu = 1 }^{T}$ the matrix valued
spectral density function of the $T$-variate  stationary
sequence $\vec \zeta(n)=(\zeta_1(n),\dots,\zeta_T(n))^{\top}$ arising from the $T$-blocking (\ref{block}) of a univariate T-PC
sequence $\zeta(n)$.

\section{Hilbert space projection method of linear interpolation}

Let $\vec {\xi}(j)$ and $\vec {\eta}(j)$ be uncorrelated T-variate
stationary stochastic sequences with the spectral density matrices
$f^{\vec \xi}(\lambda)=\left\{ f^{\vec \xi}_{\nu\mu}(\lambda)
\right\}_{\nu,\mu = 1 }^{T}$ and $f^{\vec \eta}(\lambda)=\left\{
f^{\vec \eta}_{\nu\mu}(\lambda) \right\}_{\nu,\mu=1 }^{T}$,
respectively. Consider the problem of optimal linear estimation of
the functional $$A_{s}\vec {\xi}=\sum_{l=0}^{s-1}\sum_{j=M_l+1}^{M_l+N_{l+1}}\vec {a}^{\top}(j)\vec
{\xi}(j),\,\,M_l=\sum_{k=0}^l(N_k+K_k),\,\,N_0=K_0=0,$$ that depends on the unknown values of the sequence  $\vec {\xi}(j)$,
based on observations of the sequence $\vec
{\xi}(j)+\vec{\eta}(j)$ at points $j\in {\mathbb Z}
\setminus S$, where $S=\bigcup _{l=0}^{s-1}\{M_l+1,\dots,M_l+N_{l+1}\}$.

Let the spectral densities $f^{\vec \xi}(\lambda)$ and $f^{\vec
\eta}(\lambda)$ satisfy the minimality condition
\begin{equation}\label{3.5}
\int_{-\pi}^{\pi}{Tr{\left[ {(f^{\vec \xi}(\lambda)+f^{\vec
\eta}(\lambda))^{-1}} \right]}} d\lambda <+{\infty}.
\end{equation}
\noindent Condition (\ref{3.5}) is necessary and sufficient in order that
the error-free interpolation of the unknown values of the sequence $\vec
{\xi}(j)+\vec {\eta}(j)$ is impossible \cite{Rozanov}.

 Denote by $L_{2}(f)$
the Hilbert space of vector valued functions $\vec  b(\lambda)=\left\{
b_{\nu}(\lambda) \right\}_{\nu = 1 }^{T}$ that are square integrable with
respect to a measure with the density $f(\lambda)=\left\{
f_{\nu\mu}(\lambda) \right\}_{\nu,\mu = 1 }^{T}$:
$$
\int_{-\pi}^{\pi}\vec b^{\top}(\lambda)f(\lambda)\overline{\vec b(\lambda)}d\lambda=\int_{-\pi}^{\pi}\sum_{\nu,\mu=1}^{T}b_{\nu}(\lambda) f_{\nu\mu}(\lambda) \overline{b_{\mu}(\lambda)}d\lambda<+\infty.
$$
\noindent Denote by $L_{2}^{s-}(f)$ the subspace in $L_{2}(f)$
generated by the functions $e^{ij\lambda}\delta_{\nu},\delta_{\nu}=\left\{
\delta_{\nu\mu} \right\}_{\mu = 1 }^{T}$, $\nu=1,\dots,T,\,j\in {\mathbb
Z}\setminus S$, where $\delta_{\nu\nu}=1,\delta_{\nu\mu}=0$
for $\nu\neq\mu$.

Every linear estimate $\widehat{A_s\vec {\xi}}$ of the functional
$A_{s}\vec {\xi}$ from observations of the sequence $\vec
{\xi}(j)+\vec{\eta}(j)$ at points $j\in {\mathbb Z}
\setminus S$ has the form
\begin{equation}\label{est}
\widehat{{A}_{s}\vec
{\xi}}=\int_{-\pi}^{\pi}\vec h^{\top}(e^{i\lambda})(Z^{\vec \xi}(d\lambda)+Z^{\vec \eta}(d\lambda))=\int_{-\pi}^{\pi}\sum_{\nu=1}^{T}h_{\nu}(e^{i\lambda})(Z_{\nu}^{\vec \xi}(d\lambda)+Z_{\nu}^{\vec \eta}(d\lambda)),
\end{equation}
\noindent where $Z^{\vec \xi}(\Delta)=\left\{ Z_{\nu}^{\vec \xi}(\Delta)
\right\}_{\nu = 1}^{T}$ and $Z^{\vec \eta}(\Delta)=\left\{
Z_{ \nu}^{\vec \eta}(\Delta) \right\}_{\nu = 1 }^{T}$ are orthogonal random
measures of the sequences $\vec {\xi}(j)$ and $\vec {\eta}(j)$, and
$\vec h(e^{i\lambda})=\left\{ h_{\nu}(e^{i\lambda}) \right\}_{\nu = 1
}^{T}$ is the spectral characteristic of the estimate
$\widehat{{A}_{s}\vec {\xi}}$. The function $\vec h(e^{i\lambda})\in
L_{2}^{s-}(f^{\vec \xi}+f^{\vec \eta})$.

The mean square error $\Delta(\vec h;f^{\vec \xi},f^{\vec \eta})$ of the
estimate $\widehat{{A}_{s}\vec {\xi}}$ is calculated by the formula
$$
\Delta(\vec h;f^{\vec \xi},f^{\vec \eta})=E|A_{s}\vec
{\xi}-\widehat{{A}_{s}\vec {\xi}}|^{2}=
$$
\begin{equation}\label{3.55}
=\frac{1}{2\pi}\int_{-\pi}^{\pi}{\left[ {
A_{s}(e^{i\lambda})-\vec h(e^{i\lambda})} \right]}^{\top}f^{\vec
\xi}(\lambda){\overline{\left[ {
A_{s}(e^{i\lambda})-\vec h(e^{i\lambda})} \right]}}
d\lambda+\end{equation}
$$+\frac{1}{2\pi}\int_{-\pi}^{\pi} \vec h^{\top}(e^{i\lambda})f^{\vec
\eta}(\lambda)\overline{\vec h(e^{i\lambda})}d\lambda,
$$
$$
A_{s}(e^{i\lambda})=\sum_{l=0}^{s-1}\sum_{j=M_l+1}^{M_l+N_{l+1}}\vec {a}(j)e^{ij\lambda}.
$$
\noindent The spectral characteristic $\vec h(f^{\vec \xi},f^{\vec
\eta})$ of the optimal linear estimate of $A_{s}\vec {\xi}$
minimizes the mean square error
\begin{equation}\label{3.6}
\Delta(f^{\vec \xi},f^{\vec \eta})=\Delta(\vec h(f^{\vec \xi},f^{\vec
\eta});f^{\vec \xi},f^{\vec \eta})=\mathop {\min }\limits_{\vec h \in
L_{2}^{s-}(f^{\vec \xi}+f^{\vec \eta})} \Delta (\vec h;f^{\vec
\xi},f^{\vec \eta})=\mathop {\min }\limits_{\widehat{{A}_{s}\vec
{\xi}}}E|A_{s}\vec {\xi}-\widehat{{A}_{s}\vec {\xi}}|^{2}.
\end{equation}
With the help of the
Hilbert space projection method proposed by A. N. Kolmogorov \cite{Kolmogorov} we can
find a
solution of the optimization problem (\ref{3.6}). The optimal linear estimate $\hat{A}_{s}\vec {\xi}$ is a projection of the functional $A_s \vec \xi$ on the subspace $H^-[\vec \xi+\vec \eta]=H^-[\xi_\nu(j)+\eta_\nu(j),\, j\in \mathbb{Z}\backslash S, \nu=1,\dots,T]$ of the Hilbert space $H=\{\xi: \textsf{E}\xi=0,\, \textsf{E}|\xi|^2<\infty\}$,
 generated by  values $\xi_\nu(j)+\eta_\nu(j),\, j\in \mathbb{Z}\backslash S, \nu=1,\dots,T$. The projection is characterized by the following conditions\\
 1) $ \widehat{A_s \vec \xi} \in H^-[\vec \xi+\vec \eta],$\\
 2) $A_s \vec \xi-\widehat{A_s \vec \xi}\perp H^-[\vec \xi+\vec \eta].$

Condition 2) gives us the possibility to derive the  formula for the spectral characteristic of the optimal estimate
$$
\vec h^{\top}(f^{\vec \xi},f^{\vec \eta})= \left( A^{\top}_s
(e^{i\lambda} )f^{\vec \xi}(\lambda) - C^{\top}_s (e^{i\lambda
})\right)\, \left[ f^{\vec \xi}(\lambda)+f^{\vec \eta}(\lambda)
\right]^{-1}=
$$
\begin{equation}\label{3.7}
=A_{s}^{\top}(e^{i\lambda})-\left(
{A_{s}^{\top}(e^{i\lambda})f^{\vec
\eta}(\lambda)+C_{s}^{\top}(e^{i\lambda})}\right) \left[ f^{\vec
\xi}(\lambda) +f^{\vec \eta}(\lambda)\right]^{-1},
\end{equation}
\noindent where
$$
C_{s}(e^{i\lambda})=\sum_{l=0}^{s-1} \sum_{k_l=M_l+1}^{M_l+N_{l+1}}\vec
{c}(k_l)e^{ik_l\lambda},$$
$$\vec
{c}(k_l)=\left(
c_{1}(k_l),
\dots,
c_{T}(k_l)
\right)^{\top},$$
$$l=0,\dots,s-1,\, k_l=M_l+1,\dots,M_l+N_{l+1}.
$$

Condition 1) is satisfied when the system of equalities
\begin{equation}\label{syst}
\int_{-\pi}^{\pi}\vec h^{\top}(f^{\vec \xi},f^{\vec
\eta})e^{-ij\lambda}d\lambda=0, \quad j\in S
\end{equation}
holds true.

Denote by $\textbf{D}_s$,
$\textbf{B}_s$ operators that are determined by  $T\rho\times T\rho$, $\rho=N_1+N_2+\dots+N_s$, matrices
 $$
 \textbf{D}_{s}=\begin{pmatrix}
D_{00} & D_{01} &\dots &D_{0,s-1}\\
D_{10} & D_{11} &\dots &D_{1,s-1}\\
\dots&\dots&\dots&\dots\\
D_{s-1,0} & D_{s-1,1} &\dots &D_{s-1,s-1}\\
\end{pmatrix},\,\,\,
\textbf{B}_{s}=\begin{pmatrix}
B_{00} & B_{01} &\dots &B_{0,s-1}\\
B_{10} & B_{11} &\dots &B_{1,s-1}\\
\dots&\dots&\dots&\dots\\
B_{s-1,0} & B_{s-1,1} &\dots &B_{s-1,s-1}\\
\end{pmatrix},
$$
constructed from $TN_{m+1}\times TN_{n+1}$ block-matrices $$ D_{mn}=\left\{D_{mn}(k,j)\right\}_{k=M_{m}+1}^{M_{m}+N_{m+1}} {} _{j=M_{n}+1}^{M_{n}+N_{n+1}},$$ $$B_{mn}=\left\{B_{mn}(k,j)\right\}_{k=M_{m}+1}^{M_{m}+N_{m+1}} {} _{j=M_{n}+1}^{M_{n}+N_{n+1}},$$
\[ m,n=0,\dots,s-1,\]
with elements which are
the Fourier coefficients of  the matrix functions $\left[f^{\vec \xi}(\lambda)(f^{\vec \xi}(\lambda)+f^{\vec
\eta}(\lambda))^{-1}\right]^{\top}$ and $\left[(f^{\vec
\xi}(\lambda)+f^{\vec \eta}(\lambda))^{-1}\right]^{\top}$, correspondingly:
\[
D_{mn}(k,j)=\frac{1}{2\pi}\int_{-\pi}^{\pi}\left[f^{\vec
\xi}(\lambda)(f^{\vec \xi}(\lambda)+f^{\vec
\eta}(\lambda))^{-1}\right]^{\top}e^{i(j-k)\lambda}d\lambda,\]
\[
B_{mn}(k,j)=\frac{1}{2\pi}\int_{-\pi}^{\pi}\left[(f^{\vec
\xi}(\lambda)+f^{\vec
\eta}(\lambda))^{-1}\right]^{\top}e^{i(j-k)\lambda}d\lambda,
\]
\[k=M_{m}+1
,\dots,M_{m}+N_{m+1},\]
\[j=M_{n}+1
,\dots,M_{n}+N_{n+1}.
\]
Making use of the introduced operators, relation (\ref{syst}) can be written in the  form of the equation
\[
\textbf{D}_s \vec a_s=\textbf{B}_s\vec c_s,
\]
\noindent where  $$\vec a_s=\left(\vec a^{\top}(1),\dots,\vec a^{\top}(N_1),\vec a^{\top}(M_1+1)\dots, \vec a^{\top}(M_1+N_2), \dots, \vec a^{\top}(M_{s-1}+1),\dots, \vec a^{\top}(M_{s-1}+N_s)\right)^{\top},$$
$$\vec c_s=\left(\vec c^{\top}(1),\dots,\vec c^{\top}(N_1),\vec c^{\top}(M_1+1)\dots, \vec c^{\top}(M_1+N_2), \dots, \vec c^{\top}(M_{s-1}+1),\dots, \vec c^{\top}(M_{s-1}+N_s)\right)^{\top}$$ are column-vectors.
The unknown coefficients $\vec
{c}(k_l),\,l=0,\dots,s-1,\, k_l=M_l+1,\dots,M_l+N_{l+1}$ are determined from the equation
\begin{equation}\label{c}
\vec c_s=\textbf{B}_s^{-1}\textbf{D}_s \vec a_s,
\end{equation}
where the $k_l$-th component of the vector $\vec c_s$ is calculated by the formula
$$
\vec c(k_l)=\sum _{m=0}^{s-1} \sum_{q=M_{m}+1}^{M_{m}+N_{m+1}}  C_{lm}(k_l,q) \sum _{n=0}^{s-1} \sum_{j=M_{n}+1}^{M_{n}+N_{n+1}} D_{mn}(q,j) \vec a(j),$$
$$ l=0,\dots,s-1,\,k_l=M_{l}+1
,\dots,M_{l}+N_{l+1}.
$$

The operator $\textbf{B}_s^{-1}$ is determined by $T\rho\times T\rho$ matrix
$$\textbf{B}_{s}^{-1}=\begin{pmatrix}
C_{00} & C_{01} &\dots &C_{0,s-1}\\
C_{10} & C_{11} &\dots &C_{1,s-1}\\
\dots&\dots&\dots&\dots\\
C_{s-1,0} & C_{s-1,1} &\dots &C_{s-1,s-1}\\
\end{pmatrix}
$$
that is an inverse matrix for the block-matrix $\textbf{B}_s$. Elements of $\textbf{B}_s^{-1}$ are constructed by dividing $\textbf{B}_s^{-1}$ on  $TN_{m+1}\times TN_{n+1}$ block-matrices  $C_{mn}$ and dividing each  $C_{mn}$  on $T\times T$ matrices  $C_{mn}(k,j),$ $k=M_{m}+1
,\dots,M_{m}+N_{m+1},j=M_{n}+1
,\dots,M_{n}+N_{n+1},$ $m,n=0,\dots,s-1$, in the  such way that
$$C_{mn}=\left\{C_{mn}(k,j)\right\}_{k=M_{m}+1}^{M_{m}+N_{m+1}} {} _{j=M_{n}+1}^{M_{n}+N_{n+1}}.$$

The mean-square error of the optimal estimate $\widehat{ A \vec \xi}$
is calculated by the formula (\ref{3.55}) and is of the form
\begin{equation}\label{3.8}
\Delta(f^{\vec \xi},f^{\vec
\eta})=\langle{\vec{a}_{s},\textbf{R}_{s}\vec{a}_{s}}\rangle+\langle{\vec{c}_{s},\textbf{B}_{s}\vec{c}_{s}}\rangle,
\end{equation}
\noindent where
 $\langle{a,b}\rangle$ denotes the scalar product, $\textbf{R}_{s}$ is the linear operator determined by $T\rho\times T\rho$ matrix composed with $TN_{m+1}\times TN_{n+1}$ block-matrices $$ R_{mn}=\left\{R_{mn}(k,j)\right\}_{k=M_{m}+1}^{M_{m}+N_{m+1}} {} _{j=M_{n}+1}^{M_{n}+N_{n+1}},$$
 \[m,n=0,\dots,s-1,\]
 with elements
$$
R_{mn}(k,j)=\frac{1}{2\pi}\int_{-\pi}^{\pi}{\left[ { f^{\vec
\xi}(\lambda)(f^{\vec \xi}(\lambda)+f^{\vec
\eta}(\lambda))^{-1}f^{\vec \eta}(\lambda)}
\right]}^{\top}e^{i(j-k)\lambda} d\lambda,
$$
\[k=M_{m}+1
,\dots,M_{m}+N_{m+1},\]
\[j=M_{n}+1
,\dots,M_{n}+N_{n+1}.
\]
\\
\noindent See \cite{MoklyachukMas2006a} for more details.
\\

The following statement holds true.

\begin{theorem}
\label{theorem3.1}
Let $\vec
{\xi}(j)=\left\{ \xi_{\nu}(j) \right\}_{\nu = 1 }^{T}$ and $\vec
{\eta}(j)=\left\{ \eta_{\nu}(j) \right\}_{\nu = 1}^{T}$ be
uncorrelated T-variate stationary stochastic sequences with the
spectral density matrices $f^{\vec \xi}(\lambda)=\left\{ f^{\vec
\xi}_{\nu\mu}(\lambda) \right\}_{\nu,\mu = 1 }^{T}$ and $f^{\vec
\eta}(\lambda)=\left\{ f^{\vec \eta}_{\nu\mu}(\lambda) \right\}_{\nu,\mu = 1
}^{T}$, respectively. Assume that the matrices $f^{\vec \xi}(\lambda)$
and $f^{\vec \eta}(\lambda)$ satisfy the minimality condition (\ref{3.5}).
The spectral characteristic $\vec h(f^{\vec \xi},f^{\vec \eta})$ and
the mean square error $\Delta(f^{\vec \xi},f^{\vec \eta})$ of the
optimal linear estimate of the functional $A_{s}\vec {\xi}$
based on observations of the sequence $\vec
{\xi}(j)+\vec{\eta}(j)$ at points  $j\in {\mathbb Z}
\setminus S$, are calculated by formulas (\ref{3.7}) and
(\ref{3.8}).
\end{theorem}

In the case of observations  without noise we have the following corollary.

\begin{corollary}
\label{cor3.1}
Let $\vec
{\xi}(j)=\left\{ \xi_{\nu}(j) \right\}_{\nu =1 }^{T}$ be a T-variate
stationary stochastic sequence with the spectral density matrix
$f^{\vec \xi}(\lambda)=\left\{ f^{\vec \xi}_{\nu\mu}(\lambda)
\right\}_{\nu,\mu = 1 }^{T}$, which satisfies the minimality
condition
\begin{equation}\label{3.9}
\int_{-\pi}^{\pi}{Tr{\left[ {(f^{\vec \xi}(\lambda))^{-1}} \right]}}
d\lambda <+{\infty}.\end{equation}
The spectral characteristic $\vec h(f^{\vec
\xi})$ and the mean square error $\Delta(f^{\vec \xi})$ of the
optimal linear estimate of the functional $A_{s}\vec {\xi}$
based on observations of the sequence $\vec {\xi}(j)$ at points  $j\in
{\mathbb Z} \setminus S$, are calculated by formulas
\begin{equation}\label{3.10}
\vec h^{\top}(f^{\vec
\xi})=A_{s}^{\top}(e^{i\lambda})-C_{s}^{\top}(e^{i\lambda}){\left[
{f^{\vec \xi}(\lambda)} \right]}^{-1},\end{equation}
\begin{equation}\label{3.11}
\Delta(f^{\vec
\xi})=\langle{\vec{c}_{s},\vec{a}_{s}}\rangle,\end{equation}
where $$\vec a_s=\left(\vec a^{\top}(1),\dots,\vec a^{\top}(N_1),\vec a^{\top}(M_1+1)\dots, \vec a^{\top}(M_1+N_2), \dots, \vec a^{\top}(M_{s-1}+1),\dots, \vec a^{\top}(M_{s-1}+N_s)\right)^{\top},$$
$$\vec c_s=\left(\vec c^{\top}(1),\dots,\vec c^{\top}(N_1),\vec c^{\top}(M_1+1)\dots, \vec c^{\top}(M_1+N_2), \dots, \vec c^{\top}(M_{s-1}+1),\dots, \vec c^{\top}(M_{s-1}+N_s)\right)^{\top}$$ are column-vectors and $\vec{c}_{s}=\textbf{B}_{s}^{-1}\vec{a}_{s}$. $\textbf{B}_{s}$ is a $T\rho\times T\rho$ matrix composed with $TN_{m}\times TN_n$ block-matrices $ B_{mn}=\left\{B_{mn}(k,j)\right\}_{k=M_{m-1}+1}^{M_{m-1}+N_m} {} _{j=M_{n-1}+1}^{M_{n-1}+N_n}$:
$$
B_{mn}(k,j)=\frac{1}{2\pi}\int_{-\pi}^{\pi}{\left[ { (f^{\vec
\xi}(\lambda))^{-1}} \right]}^{\top}e^{i(j-k)\lambda}
d\lambda,
$$
\[m,n=1,\dots,s,\]
\[k=M_{m-1}+1
,\dots,M_{m-1}+N_m,\]
\[j=M_{n-1}+1
,\dots,M_{n-1}+N_n.
\]
The $k_l$-th component of the vector $\vec c_s$  is calculated by the formula
$$
\vec c(k_l)=\sum _{m=0}^{s-1} \sum_{q=M_{m}+1}^{M_{m}+N_{m+1}} C_{lm}(k_l,q)  \vec a(q),$$
$$ l=0,\dots,s-1,\,k_l=M_{l}
,\dots,M_{l}+N_{l+1}.
$$
The operator $\textbf{B}_s^{-1}$ is determined by $T\rho\times T\rho$ matrix that is the inverse matrix to the block-matrix $\textbf{B}_s$. Elements of $\textbf{B}_s^{-1}$ are obtained by dividing $\textbf{B}_s^{-1}$ on  $TN_{m+1}\times TN_{n+1}$ block-matrices  $C_{mn}$ and dividing each of  $C_{mn}$  on $T\times T$ matrices  $C_{mn}(k,j),$ $k=M_{m}+1
,\dots,M_{m}+N_{m+1},j=M_{n}+1
,\dots,M_{n}+N_{n+1},$ $m,n=0,\dots,s-1$, in the  such way that
$$C_{mn}=\left\{C_{mn}(k,j)\right\}_{k=M_{m}+1}^{M_{m}+N_{m+1}} {} _{j=M_{n}+1}^{M_{n}+N_{n+1}}.$$
\end{corollary}

\begin{remark} Let $s=1,\,N_1=N$. Then
$$
A_s\vec \xi=A_N\vec \xi=\sum_{j=1}
^N \vec a^{\top}(j) \vec \xi(j).$$
The spectral characteristic $\vec h(f^{\vec
\xi},f^{\vec
\eta})$ and the mean square error $\Delta(f^{\vec \xi},f^{\vec
\eta})$ of the
optimal linear estimate of the functional $A_{N}\vec {\xi}$
based on observations of the sequence $\vec {\xi}(j)$ at points $j\in
{\mathbb Z} \setminus \{1,\dots,N\}$ with the noise  $\vec \eta (j)$  are calculated by formulas
$$
h^{\top}(f^{\vec \xi},f^{\vec \eta})= \left( A^{\top}_N
(e^{i\lambda} )f^{\vec \xi}(\lambda) - C^{\top}_N (e^{i\lambda
})\right)\, \left[ f^{\vec \xi}(\lambda)+f^{\vec \eta}(\lambda)
\right]^{-1}=
$$
$$
=A_{N}^{\top}(e^{i\lambda})-\left(
{A_{N}^{\top}(e^{i\lambda})f^{\vec
\eta}(\lambda)+C_{N}^{\top}(e^{i\lambda})}\right) \left[ f^{\vec
\xi}(\lambda) +f^{\vec \eta}(\lambda)\right]^{-1},
$$
\[
\Delta(f^{\vec \xi},f^{\vec
\eta})=\langle{\vec{a}_{N},R_{N}\vec{a}_{N}}\rangle+\langle{\vec{c}_{N},B_{N}\vec{c}_{N}}\rangle.
\]
The spectral characteristic $\vec h(f^{\vec\xi})$ and the mean square error $\Delta(f^{\vec \xi})$ of the
optimal linear estimate of the functional $A_{N}\vec {\xi}$
based on observations of the sequence $\vec {\xi}(j)$ at points $j\in
{\mathbb Z} \setminus \{1,\dots,N\}$ without noise  $\vec \eta (j)$  are calculated by formulas
$$
h^{\top}(f^{\vec
\xi})=A_{N}^{\top}(e^{i\lambda})-C_{N}^{\top}(e^{i\lambda}){\left[
{f^{\vec \xi}(\lambda)} \right]}^{-1},
$$
$$
\Delta(f^{\vec \xi})=\langle{\vec{c}_{N},\vec{a}_{N}}\rangle.
$$
For more details see \cite{Dubovetska0}, \cite{MoklyachukD2016}.
\end{remark}

\begin{example}
Let  $\vec \xi(n)=\begin{pmatrix}
\xi_{1}(n)\\
\xi_{2}(n)
\end{pmatrix}$
be a 2-variate stationary stochastic sequence.
Let $\xi_{1}(n)=\theta(n)$ be a univariate stationary sequence with
the spectral density function $f(\lambda)=\frac{1}{|1-ae^{-i\lambda}|^2},\, |a|<1$, and $\xi_{2}(n)=\theta(n)+\gamma(n)$, where $\gamma(n)$ is an uncorrelated with $\theta(n)$
univariate stationary sequence with the spectral density function
$g(\lambda)=\frac{1}{|1-be^{-i\lambda}|^{2}}, \, |b|<1$.
Consider the problem of estimation of the functional
$$
A_{2}\vec \xi=\vec \xi(1)-\vec \xi(3)=(1,1)
\begin{pmatrix}
\xi_{1}(1)\\
\xi_{2}(1)
\end{pmatrix}+(-1,-1)\begin{pmatrix}
\xi_{1}(3)\\
\xi_{2}(3)
\end{pmatrix}
$$
based on observations of $\vec \xi(n),n\in {\mathbb Z}
\setminus\{{1,3}\}$. Here $\vec a(1)=(1,1),\vec a(3)=(-1,-1)$.

In this case the spectral density matrix of $\vec{\xi}(n)$ is
$$
f^{\vec
\xi}(\lambda)=\begin{pmatrix} f(\lambda) & f(\lambda)\\
f(\lambda) & f(\lambda)+g(\lambda)
\end{pmatrix}
$$
and ${ [{f^{\vec \xi}(\lambda)}]}^{-1}$ satisfies the minimality
condition (\ref{3.9}). The matrix $\textbf{B}_{2}$ and its inverse $\textbf{B}_{2}^{-1}$,
the vector of unknown coefficients $\vec{c}_{2}$ are of the form
$$
\textbf{B}_{2}=\begin{pmatrix}
2+a^2+b^2 & -1-b^2 & 0 & 0\\
-1-b^2 & 1+b^2 & 0 & 0\\
0 & 0 & 2+a^2+b^2 & -1-b^2\\
0 & 0 & -1-b^2 & 1+b^2
\end{pmatrix},$$
$$\textbf{B}_{2}^{-1}=\begin{pmatrix}
\frac{1}{1+a^2} & \frac{1}{1+a^2} & 0 & 0\\
\frac{1}{1+a^2} & \frac{2+a^2+b^2}{(1+b^2)(1+a^2)} & 0 & 0\\
0 & 0 & \frac{1}{1+a^2} & \frac{1}{1+a^2}\\
0 & 0 & \frac{1}{1+a^2} & \frac{2+a^2+b^2}{(1+b^2)(1+a^2)}
\end{pmatrix},$$
$$\vec{c}_{2}=\left(
\frac{2}{1+a^2},
\frac{3+a^2+2b^2}{(1+b^2)(1+a^2)},
-\frac{2}{1+a^2},
-\frac{3+a^2+2b^2}{(1+b^2)(1+a^2)}
\right)^{\top}.
$$
The spectral characteristic can be calculated by (\ref{3.10})
$$
\vec h^{\top}(f^{\vec
\xi})=\left(\left(\frac{2a}{1+a^2}-\frac{b}{1+b^2}\right)+\left(\frac{b}{1+b^2}\right)e^{4i\lambda},\right.$$
$$\left.\frac{b}{1+b^2}-\frac{b}{1+b^2}e^{4i\lambda} \right).
$$
Then the optimal linear estimate of $A_{2}\vec \xi$ determined by (\ref{est}) is of the
form
$$
\widehat{A_2\vec \xi}=\left(\frac{2a}{1+a^2}-\frac{b}{1+b^2}\right)\xi_1(0)+\frac{b}{1+b^2}\xi_1(4)+\frac{b}{1+b^2}\xi_2(0)-\frac{b}{1+b^2}\xi_2(4)).
$$
The mean square error of $\hat{A}_{2}\zeta$ determined by (\ref{3.11}) is
$$\Delta(f^{\vec \xi})=\frac{8}{1+a^2}+\frac{2}{1+b^2}. $$
\end{example}

\section{Interpolation of T-PC stochastic sequences}

Let $\zeta(j)$ and $\theta(j)$ be uncorrelated T-PC stochastic
sequences. Consider the problem of optimal linear estimation  of the
functional $$A_s{\zeta}=\sum_{l=0}^{s-1}\sum_{j=M_l+1}^{M_l+N_{l+1}}{a}(j){\zeta}(j),$$
 $$
 M_l=\sum_{k=0}^l(N_k+K_k),\,\,N_0=K_0=0,$$
 that depends on the unobserved values of T-PC stochastic
sequence $\zeta(j)$, based on
observations of the sequence $\zeta(j)+\theta(j)$ at points $j\in {\mathbb
Z} \setminus S,$ where $S=\bigcup _{l=0}^{s-1}\{M_l+1,\dots,M_l+N_{l+1}\}$.

Using the  Gladyshev relation (\ref{2.2}) of PC and multivariate stationary
sequences  the problem of estimation of the functional
$A_{s}\zeta$ may be reduced to the problem of estimation of the
functional $A_{s}\vec {\xi}$ since
$$
A_{s}\zeta=\sum_{l=0}^{s-1}\sum_{j=M_l+1}^{M_l+N_{l+1}}{a}(j){\zeta}(j)=\sum_{l=0}^{s-1}\sum_{j=M_l+1}^{M_l+N_{l+1}}{a}(j)\sum_{\nu=1}^{T}e^{2\pi
ij\nu/T}\xi_{\nu}(j)=
$$
$$
=\sum_{l=0}^{s-1}\sum_{j=M_l+1}^{M_l+N_{l+1}}\sum_{\nu=1}^{T}{a}(j)e^{2\pi
ij\nu/T}\xi_{\nu}(j)=\sum_{l=0}^{s-1}\sum_{j=M_l+1}^{M_l+N_{l+1}}\vec {a}^{\top}(j)\vec
{\xi}(j)=A_{s}\vec {\xi},
$$
\noindent where $$\vec {a}^{\top}(j)=\left(
a_{1}(j),
\dots,
a_{T}(j)
\right),\,a_{\nu}(j)={a}(j)e^{2\pi
ij\nu/T},\,\nu=1,\dots,T,$$ $\vec {\xi}(j)=\left\{ \xi_{\nu}(j)
\right\}_{\nu = 1}^{T}$ is a T-variate stationary stochastic
sequence that generates the PC sequence  $\zeta(j)$.

For the interpolation problem for PC sequences the following results
hold true.

\begin{theorem}
\label{theorem3.2}
Let $\zeta(j)$ and
$\theta(j)$ be uncorrelated T-PC stochastic sequences. Then the
optimal linear estimate of the functional $A_{s}\zeta$ based on
observations of the sequence $\zeta(j)+\theta(j)$ at points $j\in {\mathbb Z}
\setminus S$, is given by the formula
\smallskip
\begin{equation}\label{est1}
\widehat{{A}_{s}\zeta}=\int_{-\pi}^{\pi}\vec h^{\top}
(f^{\vec \xi},f^{\vec \eta})
(Z^{\xi}(d\lambda)+Z^{\eta}(d\lambda))=\int_{-\pi}^{\pi}\sum_{\nu=1}^{T}h_{\nu}(f^{\vec
\xi},f^{\vec \eta})(Z_{\nu}^{\xi}(d\lambda)+Z_{\nu}^{\eta}(d\lambda)),
\end{equation}
where $\vec{\xi}(j)$ and $\vec{\eta}(j)$ are
generating sequences of the sequences $\zeta(j)$ and $\theta(j)$, correspondingly.
The spectral characteristic $\vec h(f^{\vec \xi},f^{\vec \eta})$ and the
mean square error $\Delta(f^{\vec \xi},f^{\vec \eta})$ of
$\widehat{{A}_{s}\zeta}$ are calculated by formulas (\ref{3.7}) and (\ref{3.8}), where
$\vec
{a}(j)=(a_{1}(j),\dots,a_{T}(j))^{\top},a_{\nu}(j)={a}(j)e^{2\pi
ij\nu/T},\nu=1,\dots,T$.
\end{theorem}

\begin{corollary}
\label{cor3.2}
The optimal linear estimate $ \widehat{\zeta(1)}$ of the unknown value $\zeta(1)$, based on observations of the sequence $\zeta(j)+\theta(j)$ at points  $j\in {\mathbb Z}
\setminus S$ is defined by the formula (\ref{est1}). The spectral characteristic $\vec h(f^{\vec \xi},f^{\vec \eta})$ and the
mean square error $\Delta(f^{\vec \xi},f^{\vec \eta})$  of the
optimal linear estimate $\widehat{\zeta(1)}$ are calculated by formulas (\ref{3.7}) and (\ref{3.8}), where the unknown coefficients $\vec c(k_l), \, l=0,\dots,s-1,\,k_l=M_{l}+1
,\dots,M_{l}+N_{l+1}$ are defined by formulas
$$
\vec c(k_l)=\sum _{m=0}^{s-1} \sum_{q=M_{m}+1}^{M_{m}+N_{m+1}}  C_{lm}(k_l,q) D_{00}(q,1)\vec a(1),$$
where elements $C_{lm}(k_l,q),\,l,m=0,\dots,s-1,\,k_l=M_{l}+1
,\dots,M_{l}+N_{l+1}, q=M_{m}+1,\dots,M_{m}+N_{m+1}$ are determined by the same way as in Theorem \ref{theorem3.1}.
\end{corollary}

In the case  of observations without noise we have the following corollary.

\begin{corollary}
\label{cor3.3}
 Let $\zeta(j)$ be a T-PC
stochastic sequence. Then the optimal linear estimate of the
functional $A_{s}\zeta$ based on observations of the sequence
$\zeta(j)$ at points $j\in {\mathbb Z} \setminus S$, is given
by
\smallskip
\begin{equation}\label{3.12}
\widehat{{A}_{s}\zeta}=\int_{-\pi}^{\pi}\vec h^{\top}(f^{\vec
\xi})Z^{\xi}(d\lambda)=\int_{-\pi}^{\pi}\sum_{\nu=1}^{T}h_{\nu}(f^{\vec
\xi})Z_{\nu}^{\xi}(d\lambda),
\end{equation}
where $\vec{\xi}(j)$ is generating sequence of
$\zeta(j)$. The spectral characteristic $\vec h(f^{\vec \xi})$ and the
mean square error $\Delta(f^{\vec \xi})$ of $\widehat{{A}_{s}\zeta}$ are
calculated by formulas (\ref{3.10}) and (\ref{3.11}), where  $\vec
{a}(j)=(a_{1}(j),\dots,$ $a_{T}(j))^{\top}, a_{\nu}(j)={a}(j)e^{2\pi
ij\nu/T},\nu=1,\dots,T$.
\end{corollary}

\begin{corollary}
\label{cor3.4}
The optimal linear estimate $ \widehat{\zeta(1)}$ of the unknown value $\zeta(1)$, based on observations of the sequence $\zeta(j)$ at points $j\in {\mathbb Z}
\setminus S$ is defined by the formula (\ref{3.12}). The spectral characteristic $\vec h(f^{\vec \xi})$ and the
mean square error $\Delta(f^{\vec \xi})$  of the
optimal linear estimate $ \widehat{\zeta(1)}$ are calculated by formulas (\ref{3.10}) and (\ref{3.11}), where the unknown coefficients $\vec c(k_l), \, l=0,\dots,s-1,\,k_l=M_{l}+1
,\dots,M_{l}+N_{l+1}$ are defined by formulas
$$
\vec c(k_l)= C_{00}(k_l,1)\vec a(1),$$
where elements $C_{00}(k_l,1),\, l=0,\dots,s-1,\,k_l=M_{l}+1
,\dots,M_{l}+N_{l+1}$ are determined by the same way as in Corollary \ref{cor3.1}.
\end{corollary}

\section{Interpolation of T-PC stochastic sequences with special
 sets of missed observations}

Consider the problem of optimal estimation for the linear functional
$$A_s{\zeta}=\sum_{l=0}^{s-1}\sum_{j=M_l+1}^{M_l+N_{l+1}}{a}(j){\zeta}(j),\,\,M_l=\sum_{k=0}^l(N_k+K_k),\,\,N_0=K_0=0,$$
which depends on the
unobserved values of $T$-PC sequence
${\zeta}(j)$ from observations of the sequence
${\zeta}(j)+{\theta}(j)$ at points $j\in{\mathbb
Z}\setminus S$,
where the number of missed observations at each of the intervals is a multiple of the period $T$, what means that
 $$K_1=T \cdot K_1^T, K_2=T\cdot K_2^T,\dots, K_{s-1}=T \cdot K_{s-1}^T,$$
and the number of observations at each of the intervals is a multiple of $T$ $$N_1=T \cdot N_1^T, N_2=T \cdot N_2^T,\dots, N_{s}=T \cdot N_s^T,$$
and coefficients $a(j), j\in S$ are of the form
\begin{equation}\label{aj}
a(j)=a\left(\left(j-\left[\frac{j}{T}\right]T\right)+\left[\frac{j}{T}\right]T\right)=a(\nu+\tilde{j}T)=a(\tilde{j})e^{2\pi i\tilde{j}\nu/T},\,
\end{equation}
$$
\nu=1,\dots,T,\, \tilde{j}\in \tilde{S},$$
$$
\tilde{S}=\bigcup_{l=0}^{s-1}\left\{M_l^T,\dots,M_l^T+N_{l+1}^T-1\right\},
$$
where  $\nu=T$ and $\tilde{j}=\lambda-1$, if $j=T\cdot \lambda, \,\lambda \in \mathbb{Z},$ or
$$
 a(j)=a(T\cdot \lambda)=a(T+(\lambda-1)T)=a(\lambda-1)e^{2\pi i(\lambda-1)T/T},
 $$
and $M_l=T\cdot M_l^T,\, l=0,\dots,s-1.$

Using Proposition 2.3  the linear functional $A_s \zeta$ can be written as follows
$$A_s{\zeta}=\sum_{l=0}^{s-1}\sum_{j=M_l+1}^{M_l+N_{l+1}}{a}(j){\zeta}(j)=
$$
$$=\sum_{l=0}^{s-1}\sum_{\tilde{j}={M_l^T}}^{M_l^T+N_{l+1}^T-1}\sum_{\nu=1}^{T}{a}(\nu+\tilde{j}T){\zeta}(\nu+\tilde{j}T)=\sum_{l=0}^{s-1}\sum_{\tilde{j}={M_l^T}}^{M_l^T+N_{l+1}^T-1}\sum_{\nu=1}^{T}{a}_\nu(\tilde{j}){\zeta}_\nu(\tilde{j})=
$$
\begin{equation}\label{veczeta}
=\sum_{l=0}^{s-1}\sum_{\tilde{j}={M_l^T}}^{M_l^T+N_{l+1}^T-1}{\vec a}^{\top}(\tilde{j}){\vec \zeta}(\tilde{j})=A_s{\vec\zeta},
\end{equation}
where $$\vec a(\tilde{j})=(a_{1}(\tilde{j}),\dots,a_{T}(\tilde{j}))^{\top},\,\, a_{\nu}(\tilde{j})=a(\nu+\tilde{j}T)={a}(\tilde{j})e^{2\pi
i\tilde{j}\nu/T},$$
$$
\vec \zeta(\tilde{j})=(\zeta_{1}(\tilde{j}),\dots,\zeta_{T}(\tilde{j}))^{\top},\,\, \zeta_{\nu}(\tilde{j})=\zeta(\nu+\tilde{j}T),$$
$$\nu=1,\dots,T,\, \tilde{j}\in \tilde{S}
$$
and $\vec \zeta(\tilde{j}),\,\tilde{j}\in \tilde{S}$ is $T$-variate stationary sequence, obtained by the $T$-blocking (\ref{block}) of univariate $T$-PC sequence $\zeta(j), \, j \in S$.

Let  $f^{\vec \zeta}(\lambda)$ and $f^{\vec
\theta}(\lambda)$ be the
spectral density matrices of T-variate stationary sequences $\vec \zeta(j)$
and $\vec \theta(j)$, obtained from the $T$-blocking (\ref{block}) of univariate $T$-PC sequences
$\zeta(j)$ and $\theta(j)$, respectively.
\\

Taking into account the definition of the functional $A_s{\vec\zeta}$ and Theorem \ref{theorem3.1} we can verify that the following statements hold true.

\begin{theorem}
\label{theorem3.3}
Let $\zeta(j)$ and
$\theta(j)$ be uncorrelated T-PC stochastic sequences with the
spectral density matrices $f^{\vec \zeta}(\lambda)$ and $f^{\vec
\theta}(\lambda)$ of T-variate stationary sequences $\vec \zeta(j)$
and $\vec \theta(j)$, respectively. Assume that  $f^{\vec
\zeta}(\lambda)$ and $f^{\vec \theta}(\lambda)$ satisfy the
minimality condition (\ref{3.5}). Then the optimal linear estimate of
$A_s\vec {\zeta}$ based on observations of $\vec\zeta(j)+\vec\theta(j)$ at points  $j\in
{\mathbb Z} \setminus \tilde{S}$, is given by
$$
\widehat{{A}_{s}\vec \zeta}=\int_{-\pi}^{\pi}\vec h^{\top}(f^{\vec \zeta},f^{\vec
\theta})(Z^{\vec\zeta}(d\lambda)+Z^{\vec\theta}(d\lambda))=\int_{-\pi}^{\pi}\sum_{\nu=1}^{T}h_{\nu}(f^{\vec
\zeta},f^{\vec
\theta})(Z_{\nu}^{\vec \zeta}(d\lambda)+Z_{\nu}^{\vec \theta}(d\lambda)),
$$
where $Z^{\vec \zeta}(\Delta)=\left\{ Z_{\nu}^{\vec \zeta}(\Delta)
\right\}_{\nu = 1}^{T}$ and $Z^{\vec \theta}(\Delta)=\left\{
Z_{ \nu}^{\vec \theta}(\Delta) \right\}_{\nu = 1 }^{T}$ are orthogonal random
measures of the sequences $\vec \zeta(\tilde{j})$ and $\vec \theta(\tilde{j})$. The spectral
characteristic $\vec h(f^{\vec \zeta},f^{\vec \theta})$ and the mean
square error $\Delta(f^{\vec \zeta},f^{\vec \theta})$ of $\widehat{{A}_s \vec
\zeta}$ are calculated by formulas
$$
\vec h^{\top}(f^{\vec \zeta},f^{\vec \theta})= \left( \tilde{A}^{\top}_s
(e^{i\lambda} )f^{\vec \zeta}(\lambda)
-\tilde{C}^{\top}_s (e^{i\lambda })\right)\left[ f^{\vec \zeta}(\lambda)+f^{\vec \theta}(\lambda)
\right]^{-1}
$$
\begin{equation}\label{3.13}
=\tilde{A}_{s}^{\top}(e^{i\lambda})-\left( \tilde{A}_{s}^{\top}(e^{i\lambda})f^{\vec
\theta}(\lambda)+
\tilde{C}_{s}^{\top}(e^{i\lambda})\right) \left[ f^{\vec
\zeta}(\lambda)+f^{\vec \theta}(\lambda)
\right]^{-1},
\end{equation}
\begin{equation}\label{3.14}
\Delta(f^{\vec \zeta},f^{\vec
\theta})=\langle{\vec{a}^{\zeta}_{s},\textbf{R}^\zeta_{s}\vec{a}^{ \zeta}_{s}}\rangle+\langle{\vec{c}^\zeta_{s},\textbf{B}^\zeta_{s}\vec{c}^\zeta_{s}}\rangle,
\end{equation}
where
$$\vec a^{ \zeta}_s=\left(\vec a^{\top}(0),\dots,\vec a^{\top}(N_1^T-1), \dots, \vec a^{\top}(M_{s-1}^T),\dots, \vec a^{\top}(M_{s-1}^T+N_s^T-1)\right)^{\top},$$
$$\vec c^{ \zeta}_s=\left(\vec c^{\top}(0),\dots,\vec c^{\top}(N_1^T-1), \dots, \vec c^{\top}(M_{s-1}^T),\dots, \vec c^{\top}(M_{s-1}^T+N_s^T-1)\right)^{\top},$$
$\tilde{A}_s(e^{i\lambda})=\sum_{\tilde{j}\in \tilde{S}}\vec a(\tilde{j})e^{i\tilde{j}\lambda},$ $\tilde{C}_s(e^{i\lambda})=\sum_{\tilde{j}\in \tilde{S}}\vec c(\tilde{j})e^{i\tilde{j}\lambda},$
the unknown coefficients $\vec c(\tilde{j}),\, \tilde{j}\in \tilde{S}$ are determined from the relation $$\vec{c}^\zeta_{s}=(\textbf{B}^{\zeta}_{s})^{-1}\textbf{D}^\zeta_{s}\vec{a}^{\zeta }_{s},$$
operators $\textbf{B}^\zeta_{s}, \textbf{D}^\zeta_{s}, \textbf{R}_{s}^\zeta$ are determined by $\rho\times \rho$ matrices,
constructed from $T\cdot N_{m+1}^T\times T\cdot N_{n+1}^T$ block-matrices $$B^\zeta_{mn}=\left\{B^\zeta_{mn}(k,j)\right\}_{k={M_{m}^T}}^{M_{m}^T+N_{m+1}^T-1} {} _{j=M_{n}^T}^{M_{n}^T+N_{n+1}^T-1},$$    $$ D^\zeta_{mn}=\left\{D^\zeta_{mn}(k,j)\right\}_{k={M_{m}^T}}^{M_{m}^T+N_{m+1}^T-1} {} _{j=M_{n}^T}^{M_{n}^T+N_{n+1}^T-1},$$ $$R^\zeta_{mn}=\left\{R^\zeta_{mn}(k,j)\right\}_{k={M_{m}^T}}^{M_{m}^T+N_{m+1}^T-1} {} _{j=M_{n}^T}^{M_{n}^T+N_{n+1}^T-1},\, m,n=0,\dots,s-1$$ with elements:
\[
B^\zeta_{mn}(k,j)=\frac{1}{2\pi}\int_{-\pi}^{\pi}{\left[
{( f^{\vec \zeta}(\lambda)+f^{\vec \theta}(\lambda))^{-1}}
\right]}^{\top}e^{i(j-k)\lambda} d\lambda,
\]
\[
D^\zeta_{mn}(k,j)=\frac{1}{2\pi}\int_{-\pi}^{\pi}{\left[
{ f^{\vec \zeta}(\lambda)(f^{\vec \zeta}(\lambda)+f^{\vec
\theta}(\lambda))^{-1}}
\right]}^{\top}e^{i(j-k)\lambda} d\lambda,
\]
\[
R^\zeta_{mn}(k,j)=\frac{1}{
2\pi}\int_{-\pi}^{\pi}{\left[ { f^{\vec
\zeta}(\lambda)(f^{\vec \zeta}(\lambda)+f^{\vec
\theta}(\lambda))^{-1}f^{\vec \theta}(\lambda)}
\right]}^{\top}e^{i(j-k)\lambda} d\lambda,
\]
\[m,n=0,\dots,s-1,\]
\[k={M_{m}^T}
,\dots,M_{m}^T+N_{m+1}^T-1,\]
\[j={M_{n}^T}
,\dots,M_{n}^T+N_{n+1}^T-1.
\]
\end{theorem}

In the case of observations  without noise we have the following corollary.

\begin{corollary}
\label{cor3.5}
Let $\zeta(j)$ be a
T-PC stochastic sequence with the spectral density matrix $f^{\vec
\zeta}(\lambda)$ of T-variate stationary sequence $\vec \zeta(j)$.
Assume that $f^{\vec \zeta}(\lambda)$ satisfies the minimality
condition (\ref{3.9}). Then the optimal linear estimate of $A_s{\vec \zeta}$
based on observations of $\vec \zeta(j)$ at points   $j\in {\mathbb Z}
\setminus \tilde{S}$, is given by
$$
\widehat{{A}_{s}\vec\zeta}=\int_{-\pi}^{\pi}\vec h^{\top}(f^{\vec
\zeta})Z^{\vec\zeta}(d\lambda)=\int_{-\pi}^{\pi}\sum_{\nu=1}^{T}h_{\nu}(f^{\vec
\zeta})Z_{\nu}^{\vec\zeta}(d\lambda).
$$
The spectral characteristic $\vec h(f^{\vec \zeta})$ and the
mean square error $\Delta(f^{\vec \zeta})$ of $ \widehat{{A}_s \vec\zeta}$ are
calculated by formulas
\begin{equation}\label{3.15}
\vec h^{\top}(f^{\vec \zeta})= \tilde{A}^{\top}_s (e^{i\lambda}) -
\tilde{C}^{\top}_s (e^{i\lambda }) \left[ f^{\vec
\zeta}(\lambda) \right]^{-1},
\end{equation}
\begin{equation}\label{3.16}
\Delta(f^{\vec
\zeta})=\langle{\vec{c}^\zeta_{s},\vec{a}^{\zeta}_{s}}\rangle,
\end{equation}
unknown coefficients $\vec c(\tilde{j}),\, \tilde{j}\in \tilde{S}$ are determined from the relation  $$\vec{c}^\zeta_{s}=(\textbf{B}^{\zeta}_{s})^{-1}\vec{a}_{s},$$ operator $\textbf{B}^\zeta_{s}$
is a matrix composed with $\rho\times \rho$ matrix,
constructed from $T\cdot N_{m+1}^T\times T\cdot N_{n+1}^T$ block-matrices $B^\zeta_{mn}=\left\{B^\zeta_{mn}(k,j)\right\}_{k=M_{m}^T}^{M_{m}^T+N_{m+1}^T-1} {} _{j=M_{n}}^{M_{n}+N_{n+1}-1}$ with elements:
$$
B^\zeta_{mn}(j,k)=\frac{1}{2\pi}\int_{-\pi}^{\pi}{\left[
{(f^{\vec \zeta}(\lambda))^{-1}}
\right]}^{\top}e^{i(k-j)\lambda} d\lambda,
$$
\[m,n=0,\dots,s-1,\]
\[k={M_{m}^T}
,\dots,M_{m}^T+N_{m+1}^T-1,\]
\[j={M_{n}^T}
,\dots,M_{n}^T+N_{n+1}^T-1.
\]
\end{corollary}

\begin{example}
Let $\zeta(n),\,n\in \mathbb{Z}$ be a 2-PC
stochastic sequence such that  $\zeta(2n+1)=\eta(n)$ is a univariate
stationary Ornstein-Uhlenbeck sequence with the spectral density
 $f(\lambda)=\frac{1}{ |2+e^{i\lambda}|^{2}}$ and
$\zeta(2n)=\gamma(n)$ is an uncorrelated with $\eta(n)$ univariate
stationary Ornstein-Uhlenbeck sequence with the spectral density
 $g(\lambda)=\frac{1}{|3-e^{i\lambda}|^{2}}$.

Consider the problem of estimation of the functional $$A_{1}\zeta=\zeta(1)+\zeta(2)-\zeta(3)+\zeta(4).$$

Here $S=\{1,2,3,4\}$ and $N_1=4$ is a multiple of $T=2$. Rewrite $A_{1}\zeta$ in the form (\ref{veczeta})
$$
A_{1}\zeta=e^{2\pi i1\cdot0/2}\zeta(1+0\cdot 2)+e^{2\pi  i2\cdot 0/2}\zeta(2+0\cdot2)+e^{2\pi i1\cdot1/2}\zeta(1+1\cdot2)+e^{2\pi i2\cdot2/2}\zeta(2+1\cdot2)=
$$
$$
=\vec a^{\top}(0)\vec\zeta(0)+\vec a^{\top}(1)\vec \zeta(1)=A_{1}\vec\zeta,
$$
where $\vec a(\tilde{j})=(a(\tilde{j})e^{2\pi i1\cdot\tilde{j}/2},a(\tilde{j})e^{2\pi i2\cdot\tilde{j}/2})^{\top},\,$$a(0)=1,\,a(1)=1,\,$ $\vec \zeta(\tilde{j})=(\zeta(1+\tilde{j}\cdot2),\zeta(2+\tilde{j}\cdot2))^{\top},$ $\tilde{j }\in\tilde{S}=\{0,1\}.$
In this case the spectral density matrix of $\vec{\zeta}(n)$ is of the form
$$
f^{\vec
\zeta}(\lambda)=\begin{pmatrix} f(\lambda) & 0\\
0 & g(\lambda)
\end{pmatrix}
$$
and ${ [{f^{\vec \zeta}(\lambda)}]}^{-1}$ satisfies the minimality
condition (\ref{3.9}). The matrix $\textbf{B}_{1}$ and its inverse $(\textbf{B}_{1}^{\zeta})^{-1}$,
the vector of unknown coefficients $\vec{c}_{2}$ are of the form
$$
\textbf{B}^{\zeta}_{1}=\begin{pmatrix}
5 & 0 & 2 & 0\\
0 & 10 & 0 & -3\\
2 & 0 & 5 & 0\\
0 & -3 & 0 & 10
\end{pmatrix},\,(\textbf{B}_{1}^{\zeta})^{-1}=\frac{1}{273\pi}\begin{pmatrix}
65 & 0 & -26 & 0\\
0 & 30 & 0 & 9\\
-26 & 0 & 65 & 0\\
0 & 9 & 0 & 30
\end{pmatrix},$$
$$\vec{c}_{1}^{\zeta}=\frac{1}{273}\left(
91,
39,
-91,
39
\right)^{\top}.
$$
The spectral characteristic of the optimal estimate of $A_{1}\zeta$ is of the form
$$
\vec h^{\top}(f^{\vec \zeta})=\left(
-\frac{2}{3}e^{-i\lambda}+\frac{10}{3}e^{2i\lambda} ,
\frac{3}{7}e^{-i\lambda}+\frac{3}{7}e^{2i\lambda}\right),
$$
and the optimal linear estimate of $A_{1}\vec\zeta$ is of the form
$$
\widehat{{A}_{1}\vec\zeta}=-\frac{2}{3}\zeta_1(-1)+\frac{10}{3}\zeta_1(2)+\frac{3}{7}\zeta_2(-1)+\frac{3}{7}\zeta_2(2)=
$$
$$
=-\frac{2}{3}\zeta(-1)+\frac{10}{3}\zeta(5)+\frac{3}{7}\zeta(0)+\frac{3}{7}\zeta(6).
$$
The mean square error of this estimate
$$ \Delta(f^{\vec
\zeta})=\frac{20}{21}.
$$
\end{example}

\section{Minimax (robust) method of linear interpolation}

Let $f(\lambda)$ and $g(\lambda)$ be the spectral density
matrices of  $T$-variate stationary
sequences $\vec \zeta(j)$ and $\vec \theta(j)$, obtained by $T$-blocking (\ref{block}) of $T$-PC sequences  $\zeta(j)$ and $\theta(j)$, respectively.

Formulas (\ref{3.13})--(\ref{3.16}) may be applied for finding the spectral characteristic and
the mean square error of the optimal linear estimate of the functional $A_s \vec \zeta$
only under the condition that the spectral density
matrices $f(\lambda)$ and $g(\lambda)$ are exactly known.
If the density matrices are not known exactly while a set
$D=D_{f}\times D_{g}$ of possible spectral densities is given, the
minimax (robust) approach to estimation of functionals from unknown
values of stationary sequences is reasonable. In this case we find the
estimate which minimizes the mean square error for all spectral
densities from the given set simultaneously.

\begin{definition} For a given class of
pairs of spectral densities $D=D_{f}\times D_{g}$ the spectral density
matrices $f^{0}(\lambda)\in D_{f}$, $g^{0}(\lambda)\in D_{g}$ are
called  {the least favorable} in $D$ for the optimal linear
estimation of the functional $A_{s}\vec\zeta$ if
$$
\Delta(f^0,g^0)=\Delta(\vec h(f^0,g^0);f^0,g^0)={\max_{\substack
{(f,g)\in D}}} \Delta (\vec h(f,g);f,g).
$$
\end{definition}

\begin{definition}
For a given class of
pairs of spectral densities $D=D_{f}\times D_{g}$ the spectral
characteristic $\vec h^0(\lambda)$ of the optimal linear estimate of
the functional $A_{s}\vec \zeta$ is called minimax (robust) if
$$
\vec h^0(\lambda)\in H_{D}=\mathop {\bigcap }\limits_{(f,g)\in D}
L_{2}^{s-}(f+g),$$
$$\mathop {\min_{\substack {\vec h\in H_{D}}}} {\max_
{\substack {(f,g)\in D}}} \Delta(\vec h;f,g)=\mathop {\max_{\substack
{(f,g)\in D}}} \Delta(\vec h^0;f,g).
$$
\end{definition}

Taking into consideration these definitions and the obtained relations we can
verify that the following lemmas hold true.

\begin{lemma}
\label{lem4.1}
The spectral density matrices
$f^{0}(\lambda)\in D_{f}$, $g^{0}(\lambda)\in D_{g}$, that satisfy
 condition (\ref{3.5}), are the least favorable in D for the optimal
linear  estimation of $A_{s}\vec\zeta$, if the Fourier coefficients of
the matrix functions
$$
(f^0(\lambda)+g^0(\lambda))^{-1},\quad
f^0(\lambda)(f^0(\lambda)+g^0(\lambda))^{-1},
$$
$$
f^0(\lambda)(f^0(\lambda)+g^0(\lambda))^{-1}g^0(\lambda)
$$
define matrices $\textbf{B}^{\zeta}_{s}{}^0, \textbf{D}^{\zeta}_{s}{}^0, \textbf{R}^{\zeta}_{s}{}^0$, that determine a solution of the constrained optimization problem
$$
{\max_{\substack {(f,g)\in
D}}}(\langle{\vec{a}^{\zeta}_{s},\textbf{R}^{\zeta}_{s}\vec{a}^{\zeta}_{s}}\rangle+\langle{(\textbf{B}^{\zeta}_{s})^{-1}\textbf{D}^{\zeta}_{s}\vec{a}^{\zeta}_{s},\textbf{D}^{\zeta}_{s}\vec{a}^{\zeta}_{s}}\rangle)=\langle{\vec{a}^{\zeta}_{s},\textbf{R}^{\zeta}_{s}{}^0\vec{a}^{\zeta}_{s}}\rangle+\langle{(\textbf{B}^{\zeta}_{s}{}^0)^{-1}\textbf{D}^{\zeta}_{s}{}^0\vec{a}^{\zeta}_{s},\textbf{D}^{\zeta}_{s}{}^0\vec{a}^{\zeta}_{s}}\rangle.
$$
The minimax spectral characteristic
$\vec h^0=\vec h(f^0,g^0)$ is given by (\ref{3.7}), if $\vec h(f^0,g^0)\in H_{D}$.
\end{lemma}

\begin{lemma}
\label{lem4.2}
The spectral density matrix
$f^{0}(\lambda)\in D_{f}$, that satisfies  condition (\ref{3.9}), is
the least favorable in $D_{f}$ for the optimal linear estimation
of $A_{s}\vec\zeta$ based on observations of the sequence $\vec \zeta(j)$ at points
$j\in {\mathbb Z} \setminus \tilde{S}$, if the Fourier
coefficients of the matrix function $(f^0(\lambda))^{-1}$ define the matrix
$\textbf{B}^{\zeta}_{s}{}^0$, that determine a solution of the constrained optimization problem
$$
\max_{\substack {f\in D_f}}
\langle{(\textbf{B}^{\zeta}_{s})^{-1}\vec{a}^{\zeta}_{s},\vec{a}^{\zeta}_{s}}\rangle=\langle{(\textbf{B}^{\zeta}_{s}{}^0)^{-1}\vec{a}^{\zeta}_{s},\vec{a}^{\zeta}_{s}}\rangle.
$$
The minimax spectral characteristic $\vec h^0=\vec h(f^0)$
is given by (\ref{3.10}), if $\vec h(f^0)\in H_{D}$.
\end{lemma}

The least favorable spectral densities $f^{0}(\lambda)\in D_{f}$,
$g^{0}(\lambda)\in D_{g}$ and the minimax spectral characteristic
$\vec h^0=\vec  h(f^0,g^0)$ form a saddle point of the function $\Delta(\vec h;f,g)$
on the set $H_{D}\times D$. The saddle point inequalities
$$
\Delta(\vec h^0;f,g)\leq\Delta(\vec h^0;f^0,g^0)\leq\Delta(\vec h;f^0,g^0), \quad \forall
\vec h\in H_{D}, \forall f\in D_{f}, \forall g\in D_{g}
$$
\noindent hold when $\vec h^0=\vec h(f^0,g^0)$, $\vec h(f^0,g^0)\in H_{D}$ and
$(f^0,g^0)$ is a solution of the constrained optimization problem
\begin{equation} \label{copt1}
\sup\limits_{(f,g)\in D_f\times D_g}\Delta\left(\vec h(f^0,g^0);f,g\right)=\Delta\left(\vec h(f^0,g^0);f^0,g^0\right).
\end{equation}

\noindent
The linear functional $\Delta(\vec h(f^0,g^0);f,g)$
is calculated by the formula
$$\Delta(\vec h(f^0,g^0);f,g)=\frac{1}{2\pi}\int_{-\pi}^{\pi}\left(\tilde{A}_{s}(e^{i\lambda})g^0(\lambda)+
\tilde{C}^0_{s}(e^{i\lambda})\right)^{\top}
(f^0(\lambda)+g^0(\lambda))^{-1}f(\lambda)\times$$
$$
(f^0(\lambda)+g^0(\lambda))^{-1}\overline{\left(\tilde{{A}}_{s}(e^{i\lambda})g^0(\lambda)+
\tilde{{C}}^0_{s}(e^{i\lambda})\right)}d\lambda+ \frac{1}{2\pi}\int_{-\pi}^{\pi}\left(\tilde{A}_{s}(e^{i\lambda})f^0(\lambda)-
\tilde{C}^0_{s}(e^{i\lambda})\right)^{\top}\times$$
$$(f^0(\lambda)+g^0(\lambda))^{-1}g(\lambda)(f^0(\lambda)+g^0(\lambda))^{-1}
\overline{\left(\tilde{A}_{s}(e^{i\lambda})f^0(\lambda)-
\tilde{C}^0_{s}(e^{i\lambda})\right)}d\lambda.
$$

The constrained optimization problem (\ref{copt1}) is equivalent to the unconstrained optimization problem \cite{Pshenichnyj}:
\begin{equation} \label{copt2}
\Delta_D(f,g)=-\Delta(\vec h(f^0,g^0);f,g)+\delta((f,g)\left|D_f\times D_g\right.)\rightarrow \inf,
\end{equation}
where $\delta((f,g)|D_f\times D_g)$ is the indicator function of the set $D=D_f\times D_g$.
A solution of the problem (\ref{copt2}) is characterized by the condition $0 \in \partial\Delta_D(f^0,g^0),$ where $\partial\Delta_D(f^0,g^0)$ is the subdifferential of the convex functional $\Delta_D(f,g)$ at point $(f^0,g^0)$ \cite{Rockafellar}.

The form of the functional $\Delta(\vec h(f^0,g^0);f,g)$  admits finding the derivatives and differentials of the functional in the space $L_1\times L_1$. Therefore the complexity  of the optimization problem (\ref{copt2}) is determined by the complexity of calculating of subdifferentials of the indicator functions  $\delta((f,g)|D_f\times D_g)$  of the sets $D_f\times D_g$ \cite{Ioffe}.

Taking into consideration the introduced definitions and the derived relations we can verify that the following lemma holds true.

\begin{lemma}
\label{lem4.3}
Let $(f^0,g^0)$ be a solution to the optimization problem (\ref{copt2}). The spectral densities  $f^0(\lambda)$, $g^0(\lambda)$ are the least favorable in the class $D=D_f\times D_g$ and the spectral characteristic  $\vec h^0= \vec h(f^0,g^0)$ is the minimax of the optimal linear estimate of the functional  $A_s\vec{\zeta}$  if  $\vec h(f^0,g^0) \in H_D$.
\end{lemma}

In the case of estimation of the functional based on observations without noise we have the following statement.
\begin{lemma}
\label{lem4.4}
Let $f^0(\lambda)$ satisfies the
condition (\ref{3.9}) and be a solution of the constrained optimization problem
\begin{equation}\label{4.17}
\Delta(\vec h(f^0);f)\rightarrow{sup},
f(\lambda)\in D_{f},\end{equation}
\[
\Delta(\vec h(f^0);f)=\frac{1}{2\pi}\int_{-\pi}^{\pi}\left(\tilde{C}^0_{s}(e^{i\lambda})\right)^{\top}(f^0(\lambda))^{-1}f(\lambda)(f^0(\lambda))^{-1}
\overline{\left(\tilde{C}^0_{s}(e^{i\lambda})\right)}d\lambda.
\]
Then $f^0(\lambda)$ is the least favorable
spectral density matrix for the optimal linear estimation of
$A_{s}\vec\zeta$ based on observations of the sequence $\vec \zeta(j)$ at points  $j\in
{\mathbb Z} \setminus\tilde{ S}$. The minimax spectral
characteristic $\vec h^0=\vec h(f^0)$ is given by (\ref{3.10}), if $\vec h(f^0)\in
H_{D}$.
\end{lemma}

\section{The least favorable spectral densities in $D_{0}^{-}$ }

Let $\zeta(j),\, j\in \mathbb{Z}$ be $T$-PC sequence and let $\vec\zeta(j)$ be $T$-variate stationary sequence, obtained by $T$-blocking (\ref{block}) of $T$-PC sequence  $\zeta(j)$.
Assume that the number of missed observations of the functional
$A_{s}\vec\zeta$ at each of the intervals is a multiple of the period $T$  $$K_1=T \cdot K_1^T, K_2=T\cdot K_2^T,\dots, K_{s-1}=T \cdot K_{s-1}^T$$
and the number of observations at each of the intervals is a multiple of $T$ $$N_1=T \cdot N_1^T, N_2=T \cdot N_2^T,\dots, N_{s}=T \cdot N_s^T,$$
and coefficients $a(j), j\in S$ are of the form
(\ref{aj}).

Consider the problem of minimax estimation of the functional
$A_{s}\vec\zeta$ from observations of the sequence $\vec\zeta(j)$ at points $j\in {\mathbb Z}
\setminus \tilde{S}$ without noise, under the condition that the spectral
density matrix $f(\lambda)$ of $T$-variate stationary sequence $\vec \zeta(j)$ belongs to the set
$$
D_0^{-}={\left\{ {f(\lambda)|\,\frac{1}{2\pi}\int_{-\pi}^{\pi}f^{-1}(\lambda)d\lambda=P}
\right\}},
$$

\noindent where $P=\left\{ p_{\nu\mu} \right\}_{\nu,\mu=1}^{T}$ is a given
positive definite matrix and $det P\neq0$. With the help of
Lemma \ref{lem4.4} and the method of Lagrange multipliers  we can find that a solution
$f^{0}(\lambda)$ of the constrained optimization problem (\ref{4.17}) satisfy
the following relation:
\begin{equation}\label{5.18}
{\left[
(f^{0}(\lambda))^{-1}\right]}^{\top}
\tilde{C}_{s}^0(e^{i\lambda})={\left[
(f^{0}(\lambda))^{-1}\right]}^{\top} {\vec \alpha},
\end{equation}
\noindent where $\vec \alpha=(\alpha_1,\dots,\alpha_{T})^{\top}$ is a
vector of  Lagrange multipliers,
$$\tilde{C}_{s}^0(e^{i\lambda})=\sum_{l=0}^{s-1}\sum_{\tilde{j}={M_l^T}}^{M_l^T+N_{l+1}^T-1} \vec c^{}{}^0(\tilde{j})e^{i\tilde{j}\lambda},
$$
$$\vec c^{ \zeta}_s{}^0=\left(\vec c^0(0),\dots,\vec c^0\left({N_1^T-1}\right), \dots, \vec c^0(M_{s-1}^T),\dots, \vec c^0\left(M_{s-1}^T+N_s^T-1\right)\right)^{\top},$$
unknown coefficients $\vec c^0(\tilde{j}),\, \tilde{j}\in \tilde{S}$ are determined from relation $\vec c^{\zeta}_{s}{}^0=(\textbf{B}_{s}^{0})^{-1}\vec{a}^{\zeta}_{s},$ the matrix $\textbf{B}_{s}^0$ is
constructed from the Fourier coefficients
$$
B_{s}^0(k,j)=R^{\top}(j-k)=\frac{1}{2\pi}\int_{-\pi}^{\pi}{\left[
(f^{0}(\lambda))^{-1}\right]}^{\top} e^{i(j-k)\lambda}
d\lambda,\,\,k,j\in \tilde{S}
$$
of the matrix function
${\left[ (f^{0}(\lambda))^{-1}\right]}^{\top}.$

The Fourier coefficients $R(k)=R^*(-k),k\in \tilde{S}$, found from
the equation
$$\textbf{B}_s^{0} {\vec \alpha}_s={\vec a^{\zeta}}_s,
$$
\noindent  for ${\vec \alpha}_s=(\vec \alpha,\vec 0,\dots,\vec
0)^{\top}$, satisfy relation (\ref{5.18}) and $\textbf{B}^{0}_{s}{\vec c^{\zeta}_s{}^0}=\vec{a}_{s}$. From equations above we obtain that
$$
R(k)=\begin{cases}
P(\vec a(0))^{-1}{\vec a}^{\top}(k),\,\,k\in \tilde{S},\\
0,\,\, k\in \{0,\dots,M_{s-1}^T+N_s^T-1\}\backslash \tilde{S},
\end{cases}\\
$$
\noindent  where $\left[ (\vec a(0))^{-1}\right]^{\top }\cdot \vec a
(0)=1$. The equality $R(0)=P$ follows as a consequence of the restriction on the spectral densities from the class $D_0^-$.

Let the vector-valued sequence $\vec a (k),k\in\tilde{ S}$, be such
that the matrix function $$(f^0(\lambda))^{-1}=\sum_{k=-\left(M_{s-1}^T+N_s^T-1\right)}^{M_{s-1}^T+N_s^T-1}
R(k)e^{ik\lambda}$$ is positive definite and has nonzero
determinant. Then $(f^0(\lambda))^{-1}$ can be represented
in the form \cite{Hannan1974}
$$
(f^0(\lambda))^{-1}=\left(
\sum_{k=0}^{M_{s-1}^T+N_s^T-1} Q(k) e^{-ik\lambda} \right)\cdot \left( \sum_{k=0}^{M_{s-1}^T+N_s^T-1}
Q(k) e^{-ik\lambda} \right)^*,
$$
\noindent  where $Q(k)=0_{T\times T},\, k\in \{0,\dots,M_{s-1}^T+N_s^T-1\}\setminus\tilde{ S}$. Thus $f^0(\lambda)$  is the spectral
density of the vector autoregression stochastic sequence of
order $M_{s-1}^T+N_s^T-1$ generated by the equation
\begin{equation}\label{5.19}
\sum_{k=0}^{M_{s-1}^T+N_s^T-1} Q(k) \vec \zeta (n-k)=\vec \varepsilon(n),
\end{equation}
\noindent where  $\vec \varepsilon(n)$ is a vector "white noise"
sequence. The minimax spectral characteristic $\vec h(f^0)$ is given by
\begin{equation}\label{5.20}
\vec h(f^0)=- \sum_{k=1}^{M_{s-1}^T+N_s^T-1}\overline{R(k)} (P^T)^{-1}\vec
a(0)e^{-ik\lambda}.
\end{equation}

Hence the following theorem holds true.

\begin{theorem}
\label{theorem5.1}
Let the sequence $\vec a(k)=(a_1(k),a_2(k),\dots ,a_{T}(k))^T,$ $a_\nu(k)=a(k)e^{2\pi
i\nu k/T},$ $\nu=1,\dots ,T$, which determine the linear functional
$A_s\vec\zeta$ from  observations of sequence $\vec \zeta(j)$ at points $ j\in {\mathbb Z}
\setminus\tilde{ S}$, be such that the matrix function
$$\sum_{k=-\left(M_{s-1}^T+N_s^T-1\right)}^{M_{s-1}^T+N_s^T-1} R(k)e^{ik\lambda},$$ where
$$ R(k)=R^*(-k)=\begin{cases}
P(\vec a(0))^{-1}{\vec a}^{\top}(k),\,\,k\in \tilde{S},\\
0_{T\times T},\,\, k\in \{0,\dots,M_{s-1}^T+N_s^T-1\}\backslash \tilde{S},
\end{cases}\\
$$
 \noindent is positive definite and has nonzero determinant. Then
 the least favorable in the class $D_0^-$ spectral density for the
 optimal linear estimate of $A_s\vec\zeta$ is given by the formula
\begin{equation}\label{5.21}f^{0}(\lambda)=\left( \sum_{k=-\left(M_{s-1}^T+N_s^T-1\right)}^{M_{s-1}^T+N_s^T-1}
R(k) e^{ik\lambda} \right)^{-1}.\end{equation}
\noindent The minimax spectral characteristic $\vec h(f^{0})$ is given by
(\ref{5.20}).The greatest value of the mean square error of
$\widehat{{A}_s\vec\zeta}$ is calculated by the formula
\begin{equation}\label{5.22}\Delta(f^0)=<\vec c^{\zeta}_s{}^0,\vec a^{\zeta}_s>.
\end{equation}
\end{theorem}

\begin{example}
\label{5.1}
Let $\zeta(n)$ be a 2-PC
stochastic sequence. Consider the problem of minimax estimation of the functional $$A_2\zeta=5\zeta(1)+5\zeta(2)+2\zeta(5)+2\zeta(6)$$ from observation of the sequence $\zeta(j)$ at points $j\in\mathbb{Z}\setminus \{1,2,5,6\}$ on the set $D_0^{-}$ with $P=\begin{pmatrix}
23 & 22 \\
22 & 23
\end{pmatrix}$.

Rewrite $A_2\zeta$ in the form (\ref{veczeta})
$$A_2\zeta=5\zeta(1)+5\zeta(2)+2\zeta(5)+2\zeta(6)=$$
$$
=\vec a^{\top}(0)\vec \zeta(0)+\vec a^{\top}(2)\vec \zeta(2)=A_2\vec \zeta,
$$
where $\vec a(0)=(5,5)^{\top},\, \vec a(2)=(2,2)^{\top}.$
The matrix function $$\sum_{k={-2,0,2}}
R(k)e^{ik\lambda}$$ and the representation  $$\left(
\sum_{k=0,2} Q(k) e^{-ik\lambda} \right)\cdot \left( \sum_{k=0,2}
Q(k) e^{-ik\lambda} \right)^*$$ are of the form
$$
\begin{pmatrix}
9e^{-2i\lambda}+23+9e^{2i\lambda} & 9e^{-2i\lambda}+22+9e^{2i\lambda} \\
9e^{-2i\lambda}+22+9e^{2i\lambda} & 9e^{-2i\lambda}+23+9e^{2i\lambda}
\end{pmatrix}=$$
$$=\begin{pmatrix}
2+3e^{-2i\lambda} & 1+3e^{-2i\lambda} \\
1+3e^{-2i\lambda} & 2+3e^{-2i\lambda}
\end{pmatrix}\cdot \begin{pmatrix}
2+3e^{2i\lambda} & 1+3e^{2i\lambda} \\
1+3e^{2i\lambda} & 2+3e^{2i\lambda}
\end{pmatrix}.$$

The   least favorable spectral density in the class $D_0^-$  for the
 optimal linear estimate of $A_2\vec\zeta$ by (\ref{5.21}) is of the form
$$
f^0(\lambda)=\frac{1}{45-18e^{-2i\lambda}-18e^{2i\lambda}}\begin{pmatrix}
9e^{-2i\lambda}+23+9e^{2i\lambda} & -9e^{-2i\lambda}-22-9e^{2i\lambda} \\
-9e^{-2i\lambda}-22-9e^{2i\lambda} & 9e^{-2i\lambda}+23+9e^{2i\lambda}
\end{pmatrix}.$$

The minimax spectral characteristic, calculated by (\ref{5.20}), is given by the formula
$$
\vec h(f^0)=-\begin{pmatrix}
2  \\
2
\end{pmatrix}e^{-2i\lambda}.
$$
The greatest value of the mean square error of
$\widehat{{A}_2\vec \zeta}$ takes value
$$
\Delta(f^0)=\frac{10}{9}.
$$
\end{example}

\section{The least favorable spectral densities in $D_{G}^{-}$ }

Let $\zeta(j),\, j\in \mathbb{Z}$ be $T$-PC sequence and $\vec\zeta(j)$ be $T$-variate stationary sequence, obtained by $T$-blocking (\ref{block}) of $T$-PC sequence  $\zeta(j)$.
Assume that the number of missed observations of the functional
$A_{s}\vec\zeta$ at each of the intervals is a multiple of the period $T$  $$K_1=T \cdot K_1^T, K_2=T\cdot K_2^T,\dots, K_{s-1}=T \cdot K_{s-1}^T$$
and the number of observations at each of the intervals is a multiple of $T$ $$N_1=T \cdot N_1^T, N_2=T \cdot N_2^T,\dots, N_{s}=T \cdot N_s^T,$$
and coefficients $a(j), j\in S$ are of the form
(\ref{aj}).

Consider the problem of minimax estimation of the functional
$A_{s}\vec\zeta$ from observations $\vec\zeta(j)$ at points $j\in {\mathbb Z}
\setminus\tilde{ S}$ without noise, under the condition that the spectral
density matrix $f(\lambda)$ of the  vector stationary sequence $\vec \zeta(j)$  belongs to the set
$$
D_G^{-}={\left\{ {f(\lambda)|\,\frac{1}{2\pi}\int_{-\pi}^{\pi}f^{-1}(\lambda)\cos(g\lambda)d\lambda=P(g),
\, g=0,1,\dots,G}
\right\}},
$$

\noindent where the sequence of matrices $P(g)=\left\{ p_{\nu\mu}(g) \right\}_{\nu,\mu=1}^{T}, \, P(g)=P^*(g), g=0,1,\dots,G,$ is such that the  matrix function $\sum_{g=-G}^{G}
P(g)e^{ig\lambda}$ is positive definite  and has the determinant that does not equal zero.
With the help of Lemma \ref{lem4.4} and the method of Lagrange multipliers we can find that solution
$f^{0}(\lambda)$ of the constrained optimization problem (\ref{4.17}) satisfy
the following relation:
\begin{equation}\label{5.23}
{\left[
\left(f^{0}(\lambda)\right)^{-1}\right]}^{\top}
\tilde{C}_{s}^0(e^{i\lambda})={\left[
\left(f^{0}(\lambda)\right)^{-1}\right]}^{\top} \left(\sum_{g=0}^G {\vec \alpha_g}e^{ig\lambda}\right),
\end{equation}
\noindent where $\vec \alpha_g,\, g=0,1,\dots,G$ are
Lagrange multipliers. Relation (\ref{5.23}) holds true if
$$
\sum_{\tilde{j}\in \tilde{S}} \vec{c}^0(\tilde{j}) e^{i\tilde{j}\lambda}=\sum_{g=0}^G \vec \alpha_g e^{ig\lambda}.
$$

Consider two cases: $G\geq M_{s-1}^T+N_s^T-1$ and $G<  M_{s-1}^T+N_s^T-1$.

Let $G\geq  M_{s-1}^T+N_s^T-1$. Then the Fourier coefficients of the function  $\left(f^0({\lambda})^{-1}\right)^{\top}$ determine the matrix $\mathbf{B}_s^0$ and extremum problem (\ref{4.17}) is degenerate. Let
$$
\vec \alpha_{ M_{s-1}^T+N_s^T}=\dots=\vec \alpha_G=\vec 0\,\, \mbox{and}\,\,\vec \alpha_g=0,\,g\notin \tilde{S},
$$
and $\vec \alpha_0,\dots,\vec \alpha_{ M_{s-1}^T+N_s^T-1} $ find from the equation $$\mathbf{B}_s^0 \vec \alpha_s^0=\vec a^{\vec\zeta}_s,$$ where $\vec \alpha_s^0=\left(\vec \alpha_0,\dots,\vec \alpha_{ M_{s-1}^T+N_s^T-1}\right)^{\top}$. Then the least favorable is every density $f(\lambda)\in D_G^-$ and the density
\begin{equation}\label{1g}
 f^0(\lambda) =\left(\sum_{g=-G}^{G}
P(g)e^{ig\lambda}\right)^{-1}=
\end{equation}
$$
=\left(\left(\sum_{g=0}^{G}
Q(g)e^{-ig\lambda}\right)\left(\sum_{g=0}^{G}
Q(g)e^{-ig\lambda}\right)^{*}\right)^{-1}
$$
of the vector stochastic autoregression sequence of the order $G$
\begin{equation} \label{1gg}
\sum_{g=0}^{G} Q(g) \vec \zeta(l-g)=\vec \varepsilon_l.
\end{equation}

Let $G<  M_{s-1}^T+N_s^T-1$.  Then the matrix  $\mathbf{B}_s$ is defined by the Fourier coefficients of the  function $\left(f({\lambda})^{-1}\right)^{\top}$. Among them  $P(g),
g\in \{0,\dots,G\}\cap \tilde{S},$ are known and  $P(g), g\in \tilde{S}\setminus \{0,\dots,G\},$ are unknown. The unknown $\vec
\alpha_g,\, g\in \{0,\dots,G\}\cap S$ and $P(g), g\in\tilde{ S}\setminus \{0,\dots,G\}$ we find from the equation
\begin{equation}\label{g}
\mathbf{B}_s \vec \alpha_G^0=\vec a^{\vec\zeta}_s,
\end{equation}
where ${\mathbf{\vec \alpha}}_G^0=(\vec \alpha_0,\dots,\vec \alpha_{G'},\vec
0,\dots,\vec 0)^{\top}$, $G'$ is defined from the relation $\{0,\dots,G\}\cap \tilde{S}=\{0,\dots,G'\}$. The equation (\ref{g}) can be represented  as a system of the following equations
$$
\sum_{g\in\{0,\dots,G\}\cap \tilde{S}}B_s(0,g)\vec \alpha(g)=\vec a(0),
$$
$$
\vdots
$$
$$
\sum_{g\in\{0,\dots,G\}\cap \tilde{S}}B_s(G',g)\vec \alpha(g)=\vec a(G'),
$$
$$
\vdots
$$
$$
\sum_{g\in\{0,\dots,G\}\cap \tilde{S}}B_s\left( M_{s-1}^T+N_s^T-1,g\right)\vec \alpha(g)=\vec a\left( M_{s-1}^T+N_s^T-1\right).
$$
From the first $G'$ equations we can find coefficients $\vec \alpha_0,\dots,\vec \alpha_{G'}$ and from the next equations we can find matrices $P(g), g\in \tilde{S}\setminus \{0,\dots,G\}$.

If the sequence of matrices $P(g), g\in \tilde{S},$ is such that $P(g)=P^*(g), g\in S$, the matrix function $$\sum_{g=-\left( M_{s-1}^T+N_s^T-1\right)}^{ M_{s-1}^T+N_s^T-1}
P(g)e^{ig\lambda}$$ is positive-definite and has the determinant which does not equal zero identically, then the least favorable spectral density
$ f^0(\lambda)$ is defined by the formula
\begin{equation} \label{2g}
f^0(\lambda) =\left(\sum_{g=-\left( M_{s-1}^T+N_s^T-1\right)}^{ M_{s-1}^T+N_s^T-1}
P(g)e^{ig\lambda}\right)^{-1}=
\end{equation}
$$
=\left(\left(\sum_{g=0}^{ M_{s-1}^T+N_s^T-1}
Q(g)e^{-ig\lambda}\right)\left(\sum_{g=0}^{ M_{s-1}^T+N_s^T-1}
Q(g)e^{-ig\lambda}\right)^{*}\right)^{-1}
$$
and is the density of the vector stochastic autoregression sequence of order $ M_{s-1}^T+N_s^T-1$
\begin{equation} \label{2gg}
\sum_{g=0}^{ M_{s-1}^T+N_s^T-1} Q(g) \vec\zeta(l-g)=\vec \varepsilon_l.
\end{equation}

Thus, the following theorem holds true.

\begin{theorem}
\label{theorem6.1}
The spectral density (\ref{1g}) of the vector stochastic autoregression sequence (\ref{1gg}) of order $G$, that is determined by matrices $P(g), g\in \{0,1,\dots,G\}$, is the least favorable in the class $D^{-}_G$  for the optimal estimation of the functional  $A_s
\vec\zeta$ in the case  where $G\geq  M_{s-1}^T+N_s^T-1$. If $G<  M_{s-1}^T+N_s^T-1$ and solutions $P(g),
g\in \tilde{S}\cap \{0,1,\dots,G\},$ of the equation $\mathbf{B}_s \vec \alpha_G^0=\vec a^{\vec\zeta}_s$ with coefficients $P(g),
g\in \tilde{S}\backslash \{0,1,\dots,G\},$ form a positive-definite matrix function  $\sum_{g=-\left(  M_{s-1}^T+N_s^T-1\right)}^{  M_{s-1}^T+N_s^T-1}P(g)e^{i g\lambda}$, with the determinant which does not equal zero identically, then the spectral density (\ref{2g}) of the  vector stochastic autoregression sequence  (\ref{2gg})  of order
$ M_{s-1}^T+N_s^T-1$ is the least favorable in the class $D^{-}_G$. \noindent The minimax spectral characteristic $h(f^{0})$ is calculated by the formula
(\ref{3.10}).
\end{theorem}

\section{Conclusion}

We propose formulas for calculating the mean square error and the
spectral characteristic of the optimal linear estimate of the
functional $$A_s{\zeta}=\sum_{l=0}^{s-1}\sum_{j=M_l+1}^{M_l+N_{l+1}}{a}(j){\zeta}(j),\,\,M_l=\sum_{k=0}^l(N_k+K_k),\,\,N_0=K_0=0,$$ which depends
on the unobserved values of a periodically correlated stochastic
sequence ${\zeta}(j)$. Estimates are based on observations of the sequence
${\zeta}(j)+{\theta}(j)$ at points $j\in{\mathbb Z}\setminus S$, where $S=\bigcup _{l=0}^{s-1}\{M_l+1,\dots,M_l+N_{l+1}\}$.
The sequence ${\theta}(j)$ is an uncorrelated with ${\zeta}(j)$
periodically correlated stochastic sequence. This problem is
investigated in two cases. In the first case the spectral density matrices
$f(\lambda)$ and $g(\lambda)$ of the
$T$-variate stationary sequences $\vec \zeta(n)$ and $\vec \theta(n)$, obtained by $T$-blocking of $T$-PC sequences  $\zeta(j)$ and $\theta(j)$, respectively, are suppose to be known exactly.
In this case we derived formulas for calculating the spectral characteristic and the mean-square error of the optimal estimate of the functional.
In the second case where the spectral density matrices are
not exactly known while a class $D=D_{f} \times D_{g}$ of admissible
spectral densities is given.  Formulas that determine the least favorable spectral
densities and the minimax spectral characteristic of the optimal estimate of the functional
$A_{s}\zeta$ are proposed. The problem is investigated in details for two special classes of admissible spectral densities.
Some examples of application of the obtained  results for finding optimal estimates of linear functionals and determining the least favorable spectral
densities of the optimal estimates are presented.


\begin{thebibliography}{x}


\bibitem{Abraham}
\newblock  B. Abraham,
\newblock \emph{Missing observations in time series},
\newblock Communications in Statistics--Theory and Methods, vol. 10, pp. 1643--1653, 1981.

\bibitem{Aggoun2004}
 \newblock  L. Aggoun, and R. J. Elliott,
 \newblock \emph{Measure theory and filtering: introduction and applications},
 \newblock  Cambridge University Press, 2004.

\bibitem{Arov2018}
\newblock D. Z.  Arov, and H. Dym,
\newblock \emph{Multivariate prediction, de Branges spaces, and related extension and inverse problems},
\newblock  Birkh\"auser, 2018.

\bibitem{Basawa1980}
\newblock I.V. Basawa, and B.L.S. Prakasa Rao,
\newblock \emph{Statistical inference for stochastic processes},
\newblock  London: Academic Press, 1980.

\bibitem{Bennet}
\newblock W. R. Bennett,
\newblock \emph{Statistics of regenerative digital transmission},
\newblock  Bell System Technical Journal,  vol. 37, no. 6, pp. 1501--1542, 1958.

\bibitem{Bondon1}
\newblock P. Bondon,
\newblock \emph{Influence of missing values on the prediction of a stationary time series},
\newblock Journal of Time Series Analysis, vol. 26, no. 4, pp. 519-525, 2005.

\bibitem{Bondon2}
\newblock P. Bondon,
\newblock \emph{Prediction with incomplete past of a stationary process},
\newblock Stochastic Processes and Applications. vol.98, pp. 67-76, 2002.

\bibitem{Cheng1}
\newblock R. Cheng, A.G. Miamee, and M. Pourahmadi,
\newblock \emph{Some extremal problems in $L^p(w)$},
\newblock  Proceedings of the American Mathematical Society. vol.126, pp. 2333--2340, 1998.

\bibitem{Cheng2}
\newblock R. Cheng, and M. Pourahmadi,
\newblock \emph{Prediction with incomplete past and interpolation of missing values},
\newblock Statistics \& Probability Letters. vol. 33, pp. 341--346, 1996.

\bibitem{Cohen2015}
\newblock S. Cohen, and R.J. Elliott,
\newblock \emph{Stochastic calculus and applications},
\newblock Basel: Birkhauser, 2015.


\bibitem{Daniels2008}
\newblock M. J. Daniels, and J. W. Hogan,
\newblock \emph{Missing data in longitudinal studies: strategies for Bayesian modeling and sensitivity analysis},
\newblock Boca Raton: Taylor \& Francis Group, 2008.

\bibitem{Dubovetska0}
\newblock I. I. Dubovets'ka, O.Yu. Masyutka, and M.P. Moklyachuk,
\newblock \emph{Interpolation of periodically correlated stochastic sequences},
\newblock Theory of Probability and Mathematical Statistics, vol. 84, pp. 43-56, 2012.

\bibitem{Dubovetska4}
\newblock I. I. Dubovets'ka, and M. P. Moklyachuk,
\newblock \emph{Filtration of linear functionals of periodically correlated sequences},
\newblock Theory of Probability and Mathematical Statistics, vol. 86, pp. 51-64, 2013.

\bibitem{Dubovetska7}
\newblock I. I. Dubovets'ka, and M. P. Moklyachuk,
\newblock \emph{ Minimax estimation problem for periodically correlated stochastic processes},
\newblock Journal of Mathematics and System Science, vol. 3, no. 1, pp. 26-30, 2013.

\bibitem{Dubovetska8}
\newblock I. I. Dubovets'ka, and M. P. Moklyachuk,
\newblock \emph{ Extrapolation of periodically correlated processes from observations with noise},
\newblock Theory of Probability and Mathematical Statistics, vol. 88, pp. 67-83, 2014.

\bibitem{Dubovetska9}
\newblock I. I. Dubovets'ka, and M. P. Moklyachuk,
\newblock \emph{ On minimax estimation problems for periodically correlated stochastic processes},
\newblock Contemporary Mathematics and Statistics, vol.2, no. 1, pp. 123-150, 2014.

\bibitem{Gardner1975}
\newblock W. A. Gardner, and L. E. Franks,
\newblock \emph{Characterization of cyclostationary random signal processes},
\newblock IEEE Transactions on information theory, vol. IT-21, no. 1, pp. 4-14, 1975.

\bibitem{Gardner1994}
\newblock W.A.Gardner,
\newblock \emph{Cyclostationarity in communications and signal processing},
\newblock New York: IEEE Press, 1994.

\bibitem{Gardner2006}
\newblock W.A.Gardner, A. Napolitano, L. Paura,
\newblock \emph{Cyclostationarity: Half a century of research},
\newblock Signal Processing, vol. 86, pp. 639--697, 2006.

\bibitem{Glad1961}
\newblock E. G. Gladyshev,
\newblock \emph{Periodically correlated random sequences},
\newblock  Sov. Math. Dokl., vol. 2, pp. 385--388, 1961.

\bibitem{Grenander}
\newblock U. Grenander,
\newblock \emph{A prediction problem in game theory},
\newblock Arkiv f\"or Matematik, vol. 3, pp. 371--379, 1957.


\bibitem{Grewal2015}
\newblock M. S. Grewal, and A. P. Andrews,
\newblock \emph{Kalman filtering. Theory and practice with MATLAB},
\newblock Hoboken, NJ: John Wiley \& Sons, 2015.

\bibitem{Hannan1974}
\newblock  E. J. Hannan,
\newblock \emph{Multiple time series},
\newblock Wiley Series in Probability and Mathematical Statistics. New York: John Wiley \& Sons, 1970.

\bibitem{Hurd}
\newblock H. L. Hurd, and A. Miamee,
\newblock \emph{TPeriodically correlated random sequences},
\newblock Wiley Series in Probability and Statistics,  2007.

\bibitem{Ioffe}
\newblock A. D. Ioffe, and V. M. Tihomirov,
\newblock \emph{Theory of extremal problems},
\newblock  Studies in Mathematics and its Applications, Vol. 6. Amsterdam, New York, Oxford: North-Holland Publishing Company. XII, 1979.

\bibitem{Kallianpur1980}
 \newblock G. Kallianpur,
 \newblock \emph{Stochastic filtering theory},
 \newblock New York, Heidelberg, Berlin: Springer-Verlag, 1980.

\bibitem{Kasahara}
\newblock Y. Kasahara, M. Pourahmadi, and A. Inoue,
\newblock \emph{Duals of random vectors and processes with applications to prediction problems with missing values},
\newblock Statistics \& Probability Letters, vol. 79, no. 14, pp. 1637--1646, 2009.

\bibitem{Kassam}
\newblock  S. A. Kassam, and H. V. Poor,
\newblock  \emph{Robust techniques for signal processing: A survey},
\newblock Proceedings of the IEEE, vol. 73, no. 3, pp. 433--481, 1985.

\bibitem{Kolmogorov}
\newblock A. N. Kolmogorov,
\newblock \emph{Selected works by A. N. Kolmogorov. Vol. II: Probability theory and mathematical statistics. Ed. by A. N. Shiryayev},
\newblock Mathematics and Its Applications. Soviet Series. 26. Dordrecht etc. Kluwer Academic Publishers, 1992.


\bibitem{Koroliouk1}
\newblock D. Koroliouk,
\newblock  \emph{Stationary statistical experiments and the optimal estimator for a predictable component},
\newblock Journal of Mathematical Sciences, vol. 214, no.2, pp. 220--228, 2016.

\bibitem{Koroliouk2}
\newblock D. Koroliouk, V. S. Koroliuk, E. Nicolai, P. Bisegna, L. Stella, N. Rosato,
\newblock  \emph{A statistical model of macromolecules dynamics for Fluorescence Correlation Spectroscopy data analysis},
\newblock Statistics, Optimization and Informa. Comput., vol. 4, no. 3, pp. 233-?42, 2016.

\bibitem{Koroliouk3}
\newblock D. V. Koroliouk, and V. S. Koroliuk,
\newblock  \emph{Filtration of stationary Gaussian statistical experiments},
\newblock Journal of Mathematical Sciences, vol. 229, no.1, pp. 30--35, 2018.


\bibitem{Kozak}
\newblock P. S. Kozak, and M. P. Moklyachuk,
\newblock  \emph{Estimates of functionals constructed from random sequences with periodically stationary increments},
\newblock Theory of Probability and Mathematical Statistics, vol. 97, pp. 85-98, 2018.


\bibitem{Kutoyants1998}
 \newblock Kutoyants, Yu.A.
 \newblock  \emph{Statistical inference for spatial Poisson processes},
 \newblock New York: Springer, 1998.

\bibitem{Kutoyants2004}
\newblock  Kutoyants, Yu.A.
 \newblock  \emph{Statistical inference for ergodic diffusion processes},
 \newblock   London: Springer, 2004.

\bibitem{Little2002}
\newblock Little, R. J. A.; Rubin, D. B.
\newblock  \emph{Statistical analysis with missing data},
\newblock Hoboken, NJ: Wiley, 2019.


\bibitem{Liptser2001a}
\newblock  Liptser, R.S.; Shiryaev, A.N.
 \newblock  \emph{Statistics of random processes I. General theory},
\newblock New York -- Heidelberg -- Berlin: Springer, 2001.

\bibitem{Liptser2001b}
\newblock Liptser, R.S.; Shiryaev, A.N.
\newblock  \emph{Statistics of random processes II. Applications theory},
\newblock New York -- Heidelberg -- Berlin: Springer, 2001.

\bibitem{luz1}
\newblock M. M. Luz and M. P. Moklyachuk,
\newblock \emph{Interpolation of functionals of stochastic sequences with stationary increments},
\newblock Theory of Probability and Mathematical Statistics,  vol. 87, pp. 117-133, 2013.

\bibitem{luz2}
\newblock  M. M. Luz, and M. P. Moklyachuk,
\newblock \emph{Minimax-robust filtering problem for stochastic sequences with stationary increments},
\newblock Theory of Probability and Mathematical Statistics, vol. 89, pp. 127--142, 2014.

\bibitem{luz3}
\newblock  M. M. Luz, and M. P. Moklyachuk,
\newblock \emph{Minimax-robust filtering problem for stochastic sequences with stationary increments and cointegrated sequences},
 \newblock Statistics, Optimization \&  Information Computing, vol. 2, no. 3, pp. 176--199, 2014.

\bibitem{luz4}
\newblock  M. M. Luz, and M. P. Moklyachuk,
\newblock \emph{Minimax Interpolation problem for stochastic processes with stationary increments},
 \newblock  Statistics, Optimization \& Information Computing, vol. 3, no. 1, pp. 30--41, 2015.

\bibitem{luz5}
\newblock  M. M. Luz, and M. P. Moklyachuk,
\newblock \emph{ Minimax-robust prediction problem for stochastic sequences with stationary increments and cointegrated sequences},
 \newblock Statistics, Optimization \&  Information Computing, vol. 3, no. 2, pp. 160--188, 2015.

\bibitem{luz6}
\newblock  M. M. Luz, and M. P. Moklyachuk,
\newblock \emph{ Minimax-robust filtering problem for stochastic sequences with stationary increments and cointegrated sequences},
\newblock Cogent Mathematics, vol. 3:1167811, pp. 1--21, 2016.

\bibitem{luz7}
\newblock  M. M. Luz, and M. P. Moklyachuk,
\newblock \emph{ Minimax prediction of stochastic processes with stationary increments from observations with stationary noise},
\newblock Cogent Mathematics, vol. 3:1133219, pp. 1--17, 2016.

\bibitem{luz8}
\newblock  M. M. Luz, and M. P. Moklyachuk,
\newblock \emph{Minimax interpolation of sequences with stationary increments and cointegrated sequences},
\newblock   Modern Stochastics: Theory and Applications, vol. 3, no. 1.  pp.  59--87,  2016.

\bibitem{luz9}
\newblock  M. M. Luz, and M. P. Moklyachuk,
\newblock \emph{Filtering problem for functionals of stationary sequences},
\newblock Statistics, Optimization \& Information Computing,  vol. 4, no. 1,  pp. 68--83, 2016.

\bibitem{luz10}
\newblock  M. M. Luz, and M. P. Moklyachuk,
\newblock \emph{Minimax interpolation of stochastic processes with stationary increments from observations with noise},
\newblock Theory of Probability and Mathematical Statistics, vol. 94, pp. 121--135, 2017.

\bibitem{luz11}
\newblock  M. M. Luz, and M. P. Moklyachuk,
\newblock \emph{Estimation of stochastic processes with stationary increments and cointegrated sequences},
\newblock London: ISTE Ltd, Hoboken, NJ: John Wiley \& Sons Inc., 2019.

\bibitem{Makagon1999}
\newblock A. Makagon,
\newblock \emph{Theoretical prediction of periodically correlated sequences},
\newblock Probability and Mathematical Statistics,  vol. 19, no. 2, pp. 287--322, 1999.

\bibitem{Makagon2011}
\newblock A. Makagon, A. G. Miamee, H. Salehi, and A. R. Soltani,
\newblock \emph{Stationary sequences associated with a periodically correlated sequence},
\newblock Probability and Mathematical Statistics,  vol. 31, no. 2, pp. 263--283, 2011.

\bibitem{Masyutka:Moklyachuk:Sidei1}
\newblock  O. Yu. Masyutka, M. P. Moklyachuk, and M. I. Sidei,
\newblock \emph{Interpolation problem for multidimensional stationary sequences with missing observations},
\newblock Stochastic Modeling and Applications, vol. 22, no. 2, pp. 85--103, 2018.

\bibitem{Masyutka:Moklyachuk:Sidei2}
\newblock  O. Yu. Masyutka, M. P. Moklyachuk, and M. I. Sidei,
\newblock \emph{Interpolation problem for stationary sequences with missing observations},
\newblock Statistics, Optimization \& Information Computing, vol. 7, no. 1, pp. 97-117, 2019.

\bibitem{Masyutka:Moklyachuk:Sidei3}
\newblock  O. Yu. Masyutka, M. P. Moklyachuk, and M. I. Sidei,
\newblock \emph{Interpolation problem for multidimensional stationary processes with missing observations},
\newblock Statistics, Optimization and Information Computing, vol. 7, no. 1, pp. 118-132, 2019.

\bibitem{Masyutka:Moklyachuk:Sidei4}
\newblock  O. Yu. Masyutka, M. P. Moklyachuk, and M. I. Sidei,
\newblock \emph{Filtering of multidimensional stationary processes with missing observations},
\newblock  Universal Journal of Mathematics and Applications, vol. 2, no.1, pp. 24--35, 2019.

\bibitem{Masyutka:Moklyachuk:Sidei5}
\newblock  O. Yu. Masyutka, M. P. Moklyachuk, and M. I. Sidei,
\newblock \emph{Filtering of multidimensional stationary sequences with missing observations},
\newblock Carpathian Mathematical Publications, vol.11, no.2, pp. 361-378, 2019.

\bibitem{McKnight2007}
\newblock  P. E. McKnight, K. M. McKnight, S. Sidani, and A. J. Figueredo,
\newblock \emph{Missing data: A gentle introduction},
\newblock NY: Guilford Press, 2007.

\bibitem{Moklyachuk1993}
\newblock M. P. Moklyachuk,
\newblock \emph{On a filtering problem for vector--valued sequences},
\newblock  Theory of Probability and Mathematical Statistics,  vol. 47, pp. 107--114, 1993.

\bibitem{Moklyachuk:1993UMZ}
\newblock M. P. Moklyachuk,
\newblock \emph{On minimax filtration of vector processes},
\newblock Ukrainian Mathematical Journal, vol. 45, no.3, pp. 414--423, 1993.

\bibitem{Moklyachuk:1995ROSE}
\newblock M. P. Moklyachuk,
\newblock \emph{On interpolation problem for vector--valued stochastic sequences},
\newblock Random Operators and Stochastic Equations, vol.3, no.1, pp. 63--74, 1995.


\bibitem{Moklyachuk2000}
\newblock M. P. Moklyachuk,
\newblock \emph{ Robust procedures in time series analysis},
\newblock Theory of Stochastic Processes, vol. 6, no. 3-4, pp. 127-147, 2000.

\bibitem{Moklyachuk2001}
\newblock M. P. Moklyachuk,
\newblock \emph{Game theory and convex optimization methods in robust estimation problems},
\newblock Theory of Stochastic Processes, vol. 7, no. 1-2, pp. 253--264, 2001.

\bibitem{Moklyachuk2008}
\newblock M. P. Moklyachuk,
\newblock \emph{Robust estimations of functionals of stochastic processes},
\newblock Kyiv University, Kyiv, 2008.

\bibitem{Moklyachuk2015}
\newblock M. P. Moklyachuk,
\newblock \emph{Minimax-robust estimation problems for stationary stochastic sequences},
\newblock Statistics, Optimization \& Information Computing, vol. 3, no. 4, pp. 348--419, 2015.

\bibitem{MoklyachukD2016}
\newblock M. P. Moklyachuk, and I. I. Golichenko,
\newblock \emph{Periodically correlated processes estimates},
\newblock LAP Lambert Academic Publishing, 2016.

\bibitem{MoklyachukMas2006a}
\newblock M. P. Moklyachuk, and O. Yu. Masyutka,
\newblock \emph{Interpolation of multidimensional stationary sequences},
\newblock Theory of Probability and Mathematical Statistics, vol. 73, pp. 125--133, 2006.

\bibitem{MoklyachukMas2006b}
\newblock M. P. Moklyachuk, and O. Yu. Masyutka,
\newblock \emph{Extrapolation of multidimensional stationary processes},
\newblock Random Operators and Stochastic Equations, vol. 14, pp. 233--244, 2006.

\bibitem{MoklyachukMas2008}
\newblock M. P. Moklyachuk, and O. Yu. Masyutka,
\newblock \emph{Minimax prediction problem for multidimensional stationary stochastic processes},
\newblock Communications in Statistics -- Theory and Methods, vol. 40, no. 19-20, pp. 3700--3710, 2011.

\bibitem{MoklyachukMas2012}
\newblock M. P. Moklyachuk, and O. Yu. Masyutka,
\newblock \emph{Minimax-robust estimation technique for stationary stochastic processes},
\newblock LAP Lambert Academic Publishing, 2012.

\bibitem{Moklyachuk:Mas:Gol2018}
\newblock Moklyachuk, M. P.; Masyutka, A. Yu.; Golichenko, I.I.
Estimates of periodically correlated isotropic random fields.
\newblock Nova Science Publishers Inc. New York, 2018.

\bibitem{MSidei2015}
\newblock M. P. Moklyachuk, and M. I. Sidei,
\newblock \emph{Interpolation problem for stationary sequences with missing observations},
\newblock Statistics, Optimization \& Information Computing, vol. 3, no. 3, pp. 259-275, 2015.

\bibitem{MSidei2016}
\newblock M. P. Moklyachuk, and M. I. Sidei,
\newblock \emph{Interpolation of stationary sequences observed with the noise},
\newblock Theory of Probability and Mathematical Statistics, vol. 93, pp. 143-156, 2016.

\bibitem{MSidei2016}
\newblock M. P. Moklyachuk, and M. I. Sidei,
\newblock \emph{Filtering problem for stationary sequences with missing observations}.
\newblock Statistics, Optimization \& Information Computing,  vol. 4, no. 4, pp. 308 - 325, 2016.

\bibitem{MSidei2016b}
\newblock M. Moklyachuk, and M. Sidei,
\newblock \emph{Filtering Problem for functionals of stationary processes with missing observations},
\newblock Communications in Optimization Theory, 2016, pp.1-18, Article ID 21, 2016.

\bibitem{MSidei2017}
\newblock M. P. Moklyachuk, and M. I. Sidei,
\newblock \emph{Extrapolation problem for stationary sequences with missing observations},
\newblock Statistics, Optimization \& Information Computing,  vol. 5, no. 3, pp. 212--233, 2017.

\bibitem{Moklyachuk:Mas:Sidei}
\newblock  M. P. Moklyachuk,  M. I. Sidei, and O. Yu. Masyutka,
\newblock \emph{Estimation of stochastic processes with missing observations},
\newblock New York: Nova Science Publishers, 2019.

\bibitem{Nap1}
\newblock A. Napolitano,
\newblock \emph{Cyclostationarity: Limits and generalizations},
\newblock Signal processing, vol. 120, pp. 323--347, 2016.

\bibitem{Nap2}
\newblock A. Napolitano,
\newblock \emph{Cyclostationarity: New trends and applications},
\newblock Signal processing, vol. 120, pp. 385--408, 2016.

\bibitem{Pelagatti}
\newblock  M. M. Pelagatti,
\newblock \emph{Time series modelling with unobserved components}
\newblock New York: CRC Press, 2015.

\bibitem{PrakasaRao2010}
 \newblock B.L.S. Prakasa Rao,
 \newblock \emph{Statistical inference for fractional diffusion processes},
 \newblock Chichester: John Wiley \& Sons,  2010.

\bibitem{PrakasaRao2012}
 \newblock B.L.S. Prakasa Rao,
 \newblock \emph{Associated sequences, demimartingales and nonparametric inference}
 \newblock Basel: Birkh\"auser, 2012.

\bibitem{Pourahmadi}
\newblock M. Pourahmadi, A. Inoue, and Y. Kasahara
\newblock \emph{A prediction problem in $L^2(w)$}.
\newblock Proceedings of the American Mathematical Society. Vol. 135, No. 4, pp. 1233-1239, 2007.

\bibitem{Pshenichnyj}
\newblock B. N. Pshenichnyj,
\newblock \emph{Necessary conditions of an extremum},
\newblock Pure and Applied mathematics. 4. New York: Marcel Dekker, 1971.

\bibitem{Rajarshi2012}
\newblock M. B. Rajarshi,
\newblock \emph{Statistical inference for discrete time stochastic processes},
\newblock New Dehli: Springer India,  2012.

\bibitem{Rockafellar}
\newblock R. T. Rockafellar,
\newblock \emph{Convex Analysis},
\newblock Princeton University Press, 1997.

\bibitem{Rozanov}
\newblock Yu. A. Rozanov,
\newblock \emph{Stationary stochastic processes},
\newblock San Francisco-Cambridge-London-Amsterdam: Holden-Day, 1967.

\bibitem{Rozovsky2018}
\newblock B. L. Rozovsky, S. V. Lototsky,
\newblock \emph{ Stochastic evolution systems. Linear theory and applications to non-linear filtering},
\newblock  Springer, 2018.

\bibitem{VastPoor}
\newblock  K. S. Vastola, and H. V. Poor,
\newblock  \emph{ An analysis of the effects of spectral uncertainty on Wiener filtering},
\newblock  Automatica, vol. 28, pp. 289--293, 1983.

\bibitem{Wiener}
\newblock  N. Wiener,
\newblock \emph{ Extrapolation, interpolation and smoothing of stationary time series. With engineering applications},
\newblock  The M. I. T. Press, Massachusetts Institute of Technology, Cambridge, Mass., 1966.

\bibitem{Woodward2017}
\newblock  Woodward, W.A.; Gray, H.L.; Elliott, A.C.
\newblock \emph{Applied time series analysis with R},
\newblock  CRC Press, Taylor \& Francis Group, 2017.


\bibitem{Yaglom-a}
\newblock A. M. Yaglom,
\newblock \emph{Correlation theory of stationary and
related random functions. Vol. 1: Basic results},
\newblock Springer Series in Statistics, Springer-Verlag, New York etc., 1987.

\bibitem{Yaglom-b}
\newblock A. M. Yaglom,
\newblock \emph{ Correlation theory of stationary and
related random functions. Vol. 2: Supplementary notes and
references},
\newblock Springer Series in Statistics, Springer-Verlag, New York etc., 1987.



\end{thebibliography}
 \end{document}